\documentclass[final,3p]{elsarticle}
\makeatletter 
\def\ps@pprintTitle{%
 \def\@oddfoot{\footnotesize\itshape
       To be published in \ifx\@journal\@empty Elsevier
       \else\@journal\fi\hfill \today}}
\makeatother

\usepackage{amsmath,amsfonts,amssymb} 
\usepackage{mathtools} 
\usepackage{upgreek} 
\usepackage{easybmat} 
\usepackage{fixmath} 
\usepackage[font=footnotesize,list=true]{subcaption} 
\usepackage{booktabs} 
\usepackage[exponent-product=\cdot]{siunitx} 
\renewcommand{\ang}[1]{{{#1}^\circ}}
\usepackage{bm} 
\usepackage{xcolor}
\usepackage{microtype}
\usepackage{lmodern}

\renewcommand{\vec}[1]{\mathbold #1} 
\newcommand{\besselJ}{\mathrm{J}} 	

\newcommand{\idiff}{\, \mathrm{d}} 
\newcommand{\zerovec}{\mathbf{0}}
\newcommand{\transpose}{\top}

\newcommand{\R}{\mathbb{R}}
\providecommand{\C}{\mathbb{C}}

\newcommand{\N}{\mathbb{N}}
\newcommand{\PI}{\uppi}
\newcommand{\euler}{\mathrm{e}}
\newcommand{\imag}{\mathrm{i}}

\providecommand*{\pderiv}[3][]{\frac{\partial^{#1}#2}{\partial #3^{#1}}}

\providecommand*{\ppderiv}[3]{\frac{\partial^2{#1}}{\partial #2 \partial #3}}
\providecommand*{\pppderiv}[4]{\frac{\partial^3{#1}}{\partial #2 \partial #3 \partial #4}}
\renewcommand{\geq}{\geqslant}
\renewcommand{\leq}{\leqslant}

\DeclareMathSymbol{\GAMMA}{\mathalpha}{operators}{0}

\DeclareMathOperator{\TS}{TS}

\renewcommand\Re{\operatorname{Re}}

\newcommand{\diag}[1]{\operatorname{diag}\left(#1\right)}

\renewcommand{\Xi}{{\vec{t}_1}}

\newcommand{\energyNorm}[2]{%
  {\left\vert\kern-0.25ex\left\vert\kern-0.25ex\left\vert #1 
    \right\vert\kern-0.25ex\right\vert\kern-0.25ex\right\vert}_{#2}
}

\let\originalleft\left
\let\originalright\right
\renewcommand{\left}{\mathopen{}\mathclose\bgroup\originalleft}
\renewcommand{\right}{\aftergroup\egroup\originalright}

\newcommand{\st}{\,:\,}

\usepackage{xspace}

\usepackage[hidelinks,pdftex]{hyperref} 
\hypersetup{colorlinks = true, citecolor=blue, linkcolor=blue, urlcolor=blue} 
\usepackage{cleveref} 
\crefname{subsection}{subsection}{subsections}
\crefname{equation}{Eq.}{Eqs.}

\usepackage{longtable}
\usepackage{enumitem}

\journal{Computer methods in applied mathematics and engineering}

\begin{document}

\begin{frontmatter}

\title{Isogeometric Analysis of Acoustic Scattering with Perfectly Matched Layers (IGAPML)}

\author[SINTEF]{Jon Vegard Ven{\aa}s\texorpdfstring{\corref{cor}}{}}
\ead{JonVegard.Venas@sintef.no}

\author[NTNU,SINTEF]{Trond Kvamsdal}
\ead{Trond.Kvamsdal@sintef.no}

\address[SINTEF]{SINTEF Digital, Mathematics and Cybernetics, 7037 Trondheim, Norway}
\address[NTNU]{Department of Mathematical Sciences, Norwegian University of Science and Technology, 7034 Trondheim, Norway}

\cortext[cor]{Corresponding author.}

\begin{abstract}
The perfectly matched layer (PML) formulation is a prominent way of handling radiation problems in unbounded domain and has gained interest due to its simple implementation in finite element codes. However, its simplicity can be advanced further using the isogeometric framework. This work presents a spline based PML formulation which avoids additional coordinate transformation as the formulation is based on the same space in which the numerical solution is sought. The procedure can be automated for any convex artificial boundary. This removes restrictions on the domain construction using PML and can therefore reduce computational cost and improve mesh quality. The usage of spline basis functions with higher continuity also improves the accuracy of the PML-approximation and the numerical solution.
\end{abstract}

\end{frontmatter}

\section{Introduction}
\begin{figure}
	\centering
	\includegraphics[scale=1]{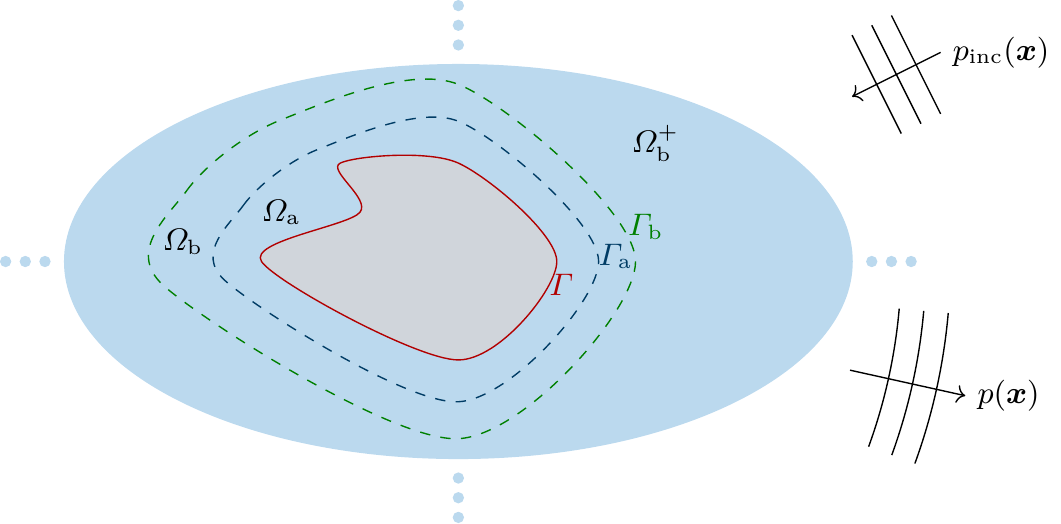}
	\caption[Illustration of artificial boundary]{Two convex boundaries $\Gamma_{\mathrm{a}}$ and $\Gamma_{\mathrm{b}}$ defines the PML around the scatterer defined by $\Gamma$ such that the exterior domain $\Omega^+$ is decomposed by the three domains $\Omega_{\mathrm{a}}$ (which is bounded by $\Gamma$ and $\Gamma_{\mathrm{a}}$), $\Omega_{\mathrm{b}}$ (which is bounded by $\Gamma_{\mathrm{a}}$ and $\Gamma_{\mathrm{b}}$) and $\Omega_{\mathrm{b}}^+$. Thus, $\Omega^+ = \Omega_{\mathrm{a}} \cup \Omega_{\mathrm{b}} \cup \Omega_{\mathrm{b}}^+$.}
	\label{Fig:artificialBoundary}
\end{figure}

Scattering problems involve \textit{unbounded exterior domains}, $\Omega^+$ (see \Cref{Fig:artificialBoundary}). 
Boundary Element Method (BEM) is a popular approach for solving such problems~\cite{Sauter2011bem,Schanz2007bea,Marburg2008cao,Chandler_Wilde2012nab}. Alternatively, a common approach for solving such problems with the finite element method (FEM) is to introduce an artificial boundary that encloses the scatterer. On the artificial boundary some sort of absorbing boundary condition (ABC) is prescribed. The problem is then reduced to a finite domain, the bounded domain between the scatterer and the artificial boundary can then be discretized with finite elements. Several methods exist for handling the exterior Helmholtz problem (on unbounded domain), including a) the perfectly matched layer (PML) method after B{\'e}renger~\cite{Berenger1994apm,Berenger1996pml,Matuszyk2012pfe}, b) Dirichlet to Neumann-operators (DtN-operators)~\cite{Givoli2013nmf}, c) local differential ABC operators~\cite{Shirron1995soe,Bayliss1982bcf,Hagstrom1998afo,Tezaur2001tdf}, and d) the infinite element method (IEM)~\cite{Bettess1977ie,Bettess1977dar,Demkowicz2001aoa}.

In earlier works we have developed isogeometric (IGA) methods~\cite{Hughes2005iac} for the IEM~\cite{Venas2018iao} and the BEM~\cite{Venas2020ibe} approaches and achieved significant improved accuracy compared to use of $C^0$ continuous FEM due to the increased inter-element continuity of the splines basis functions.

Regarding IGA for acoustic scatterings most authors have developed methods for BEM. Simpson and coworkers coined the word  {\em IGABEM\/} for isogeometric methods for BEM in~\cite{Simpson2012atd} and presented their first paper on IGABEM for acoustic scattering two years later~\cite{Simpson2014aib}. However, the first paper on isogeometric BEM~\cite{Peake2013eib} was published a year before by M. J. Peake during his PhD-study at University of Durham (UK) under supervision of Prof. J. Trevelyn and Prof. G. Coates. Here, the so-called eXtended Isogeometric Boundary Element Method (XIBEM) was introduced and further developed in the two follow up papers~\cite{Peake2014eai,Peake2015eib}.
Inspired by these initial papers several investigations of IGABEM applied to acoustics have been pursued by different authors~\cite{Dolz2016aib,Coox2017aii,Keuchel2017eoh,Dolz2018afi,Sun2019dib,Venas2020ibe,Wu2020iib}.
Recently, acoustic optimization and IGABEM has been pursued with success by a few groups, see e.g.,~\cite{Liu2017soo,Chen2018aia,Chen2019ifm,Chen2019sso,Chen2020ato, Shaaban2020sob,Shaaban2020ibe,Shaaban20213ib,Shaaban2022aib}. 

The boundary element method (BEM) avoids introducing an artificial boundary as it only relies on a computational domain on the surface of the scatterer. Moreover, solid domains are usually represented by surfaces in CAD-systems, such that if modeling of an elastic scatterer using IGA with the same spline basis as the CAD-model, the BEM does not need a surface-to-volume parametrization. This represents a significant advantage compared to the other approaches regarding interoperability between  design and analysis. Thus, the popularity among the IGA community to develop isogeometric methods for BEM (IGABEM) is understandable.

However, we experienced significant challenges related to numerical integration, fictitious eigenfrequencies, and memory requirements and solution times of the resulting algebraic system. These topics are current research areas, and we refer to~\cite{Xie2021aam} for references. Furthermore, the frequency spectra of excitation generally have a broad frequency band and thus multi-frequency analysis is often required. Because of the frequency-dependent property, both the traditional BEM and fast multipole BEM must be applied to recalculate all the entries in the system. One way to circumvent the difficulty incurred by acoustic frequency sweeps is to use reduced order models (ROM)~\cite{Hetmaniuk2012raa,Hetmaniuk2013aas,Quarteroni2015rbm}, but for BEM we must overcome the following two challenges (1) how to construct an orthonormal basis and (2) how to avoid the assembly of system matrices for each frequency before projection. Reduced order modeling of BEM for acoustic scattering is recently addressed in~\cite{Xie2021aam} but are still not yet a matured computational methodology.

Thus, use of IGABEM for addressing acoustic scattering is a versatile but challenging computational methodology. In particular, for efficient handling of frequency sweeps by means of reduced order models (ROM) it seems to be of interest to investigate alternative classical isogeometric finite element methods.

The IEM is very efficient for cases where we can locate the artificial boundary close to the scatterer and represent it with ellipsoidal coordinate systems. We developed isogeometric methods for IEM (hereafter denoted {\em IGAIEM\/} in~\cite{Venas2018iao}\footnote{To the best of our knowledge, our paper seems to be the only one combining isogeometric methods with infinite elements.} and achieved significant improved accuracy compared to use of $C^0$ continuous FEM. However, in general the need for a surface-to-volume parametrization between the scatterer and the artificial boundary is a disadvantage. Furthermore, we experienced severe challenges with high condition numbers of the system matrix when the number of radial shape functions in the infinite elements is large. This becomes a problem for more complex geometries as the number of radial shape functions must be increased to achieve higher precision. Again, there might be remedies for reducing the conditioning number, see e.g.~\cite{Safjan2002tdi} where this have been addressed for $C^0$-Lagrange FE, and that is something we will address in an upcoming paper on IGAIEM.

Compared to IEM the PML approach is not prone to ill-conditioning of the system matrix. For most applications the accuracy of the PML is comparable to IEM (for engineering precision; 1\% relative error in energy norm). 

Unfortunately, the effective implementation of the PML for convex domains of general shape has not been straight forward because of the geometric parameters that has been required to define the PML-domain. However, B\'{e}riot and Modave~\cite{Beriot2020aap} have recently presented a method for $C^0$-Lagrange finite elements that simplifies the implementations significantly. It builds upon the idea of locally conformal PML layers~\cite{Ozgun2007nml,Ozgun2007nfp}. In the present work we will investigate the use of IGA inspired by this idea to develop what we denote as {\em IGAPML\/}. Thus, the IGAPML developed herein enables us to choose the artificial boundary to be an arbitrary convex boundary represented by a NURBS parametrization. That is, we are not restricted to domains defined by the ellipsoidal, cylindrical, or Cartesian coordinate system. 

In recent papers~\cite{Mi2021ilc,Drzisga2020tsm} PML in the isogeometric framework has been presented. The present approach enables a generalization to the NURBS parametrization for the PML layer and considers different stretching functions. The stretching function recommended in the present work not only gives improved results but also reduces the number of PML parameters to tune.

A challenge for the present approach is the required surface-to-volume parametrizations from boundary representations of complex industrial CAD-models. This problem might contain trimmed NURBS patches and non-watertight models subject to a CAD cleanup, see~\cite{Marussig2017aro} for a comprehensive review. However, this challenge is ongoing research (e.g.~\cite{Urick2019wbo,Hiemstra2020tun}) and is considered out of scope for this article. The present work has focused on the automated construction of the PML-layer \textit{given} such a volumetric parametrization.
\section{Perfectly matched layer (PML) for exterior Helmholtz problems}
\label{Sec2:exteriorHelmholtz}
We partition the unbounded domain $\Omega^+$ into three domains by the boundaries $\Gamma_{\mathrm{a}}$ and $\Gamma_{\mathrm{b}}$; $\Omega_{\mathrm{a}}$, $\Omega_{\mathrm{b}}$ and $\Omega_{\mathrm{b}}^+$, see \Cref{Fig:artificialBoundary}. Due to the absorbing property of the PML layer, $\Omega_{\mathrm{b}}$, only $\Omega_{\mathrm{a}}$ and $\Omega_{\mathrm{b}}$ need to be discretized by finite elements.

The exterior Helmholtz problem is given by (with wavenumber $k$)
\begin{alignat}{3}
	\nabla^2 p + k^2 p &= 0 	&&\text{in}\quad \Omega^+,\label{Eq:HelmholtzEqn}\\
	\partial_n p &= g						&&\text{on}\quad \Gamma,\label{Eq:HelmholtzEqnNeumannCond}\\
	\pderiv{p}{r}-\imag k p &= o\left(r^{-1}\right)\quad &&\text{with}\quad r=|\vec{x}|\label{Eq:sommerfeldCond}
\end{alignat}
where the Sommerfeld condition~\cite{Sommerfeld1949pde} in \Cref{Eq:sommerfeldCond} restricts the field in the limit $r\to\infty$ uniformly in $\hat{\vec{x}}=\frac{\vec{x}}{r}$, such that no scattered waves, $p$, originate from infinity. The Neumann condition given by the function $g$ will in the case of rigid scattering be given by the incident wave $p_{\mathrm{inc}}$. Zero displacement of the fluid normal on the scatterer (rigid scattering) implies that $\partial_n(p+p_{\mathrm{inc}}) = 0$ where $\partial_n$ denotes the partial derivative in the normal direction on the surface $\Gamma$ (pointing ``out'' from $\Omega^+$), which implies that
\begin{equation}
	g = -\pderiv{p_{\mathrm{inc}}}{n}.
\end{equation}
In this work we consider plane incident waves (with amplitude $P_{\mathrm{inc}}$) traveling in the direction $\vec{d}_{\mathrm{s}}$, which can be written as
\begin{equation}\label{Eq:p_inc}
	p_{\mathrm{inc}} = P_{\mathrm{inc}}\euler^{\imag k\vec{d}_{\mathrm{s}}\cdot\vec{x}}.
\end{equation}

\subsection{Far field pattern}
The quantity of interest is the \textit{target strength} defined by
\begin{equation}\label{Eq:TS}
	\TS = 20\log_{10}\left(\frac{|p_0(\hat{\vec{x}})|}{|P_{\mathrm{inc}}|}\right)
\end{equation}
where the far field pattern of the scattered pressure, $p$, is given by
\begin{equation}\label{Eq:farfield}
	p_0(\hat{\vec{x}}) =  \lim_{r\to\infty} r \euler^{-\imag k r}p(r\hat{\vec{x}}),
\end{equation}
with $r = |\vec{x}|$ and $\hat{\vec{x}} = \vec{x}/|\vec{x}|$ being the far field observation point. The observation point can be represented in terms of the aspect angle $\alpha$ and elevation angle $\beta$
\begin{equation*}
	\hat{\vec{x}} = \begin{bmatrix}
		\cos\beta\cos\alpha\\
		\cos\beta\sin\alpha\\
		\sin\beta
	\end{bmatrix}.
\end{equation*}
We also use this convention in describing the direction of the incident wave
\begin{equation*}
	\vec{d}_{\mathrm{s}} = \begin{bmatrix}
		\cos\beta_{\mathrm{s}}\cos\alpha_{\mathrm{s}}\\
		\cos\beta_{\mathrm{s}}\sin\alpha_{\mathrm{s}}\\
		\sin\beta_{\mathrm{s}}
	\end{bmatrix}.
\end{equation*}
The far field pattern can be computed by (cf.~\cite[p. 32]{Ihlenburg1998fea})
\begin{equation}\label{Eq:HelmholtzIntegralFarField}
	p_0(\hat{\vec{x}}) = -\frac{1}{4\PI}\int_{\Gamma}\left[ \imag k p(\vec{y})\hat{\vec{x}}\cdot\vec{n}(\vec{y}) + \pderiv{p(\vec{y})}{n(\vec{y})}\right]\euler^{-\imag k \hat{\vec{x}}\cdot\vec{y}}\idiff \Gamma(\vec{y}).
\end{equation}
from which the target strength in~\Cref{Eq:TS} may be computed.

\subsection{Weak formulation for the Helmholtz equation}
The weak formulation is given by (the involved spaces are described in~\cite{Ihlenburg1998fea})
\begin{equation}
	\text{Find} \quad p\in H_w^{1+}(\Omega^+)\quad\text{such that}\quad B(q,p) = L(q),\qquad \forall q\in H_{w^*}^1(\Omega^+),
\end{equation}
where the bilinear form is given by
\begin{equation*}
	B(q,p) = \int_{\Omega^+} \left[\nabla q\cdot\nabla p-  k^2qp\right]\idiff\Omega
\end{equation*}
and the corresponding linear form is given by
\begin{equation*}
	L(q) = \int_\Gamma qg\idiff\Gamma.
\end{equation*}

\subsection{Truly perfectly matched layers}
\label{se:PML}
We will here develop a general spline-based (GSB) PML method which avoids intermediate transformation to spherical/cylindrical/Cartesian coordinates. The idea is to construct the spline patches such that the directions we want to have a decaying property is represented by parametric directions of the spline patch. As NURBS can represent spherical and cylindrical patches in addition to the trivial Cartesian patches, it can resolve the standard behavior obtained by the classical PML-formulations using the appurtenant coordinate systems. 

The NURBS basis is constructed using B-splines. Therefore, an understanding of B-splines is crucial to understanding NURBS~\cite{Cottrell2006iao}. 
We extend the classical~\cite{Cottrell2006iao} definition (using the Cox-de Boor formula) to evaluations in the complex parametric space as follows. 
Let $\check{p}$ be the polynomial order\footnote{The usage of a check sign above the polynomial order $p$ is to avoid ambiguity between the polynomial order and the scattered pressure.}, let $n$ be the number of basis functions and define a \textit{knot vector} $\vec{t} = \{\xi_1,\xi_2,\dots,\xi_{n+\check{p}+1}\}$ to be an ordered vector with non-decreasing elements, called \textit{knots}. Then, the $n$ B-splines, $\left\{B_{i,\check{p},\vec{t}}\right\}_{i\in [1,n]}$, are recursively defined by
\begin{equation*}
	B_{i,\check{p},\vec{t}}(\xi) = \frac{\xi-\xi_i}{\xi_{i+\check{p}}-\xi_i}B_{i,\check{p}-1,\vec{t}_1}(\xi)+\frac{\xi_{i+\check{p}+1} - \xi}{\xi_{i+\check{p}+1} -\xi_{i+1}}B_{i+1,\check{p}-1,\vec{t}_1}(\xi)
\end{equation*}
starting with 
(the only alteration to the classical Cox-de Boor formula is that we here take the real part of the parameter $\xi$)
\begin{equation}\label{Eq:orderOneBspline}
	B_{i,0,\vec{t}_1}(\xi)=\begin{cases}
		1 & \text{if } \xi_i\leq \Re\xi < \xi_{i+1}\\
		0 & \text{otherwise.}
	\end{cases}
\end{equation}
With B-splines in our arsenal, we are ready to present Non-Uniform Rational B-Splines (NURBS). Let $\{w_i\}_{i\in [1,n]}$ be a set of \textit{weights}, and define the \textit{weighting function} by
\begin{equation*}
	W(\xi) = \sum_{\tilde{i}=1}^{n} B_{\tilde{i},\check{p},\vec{t}}(\xi) w_{\tilde{i}}.
\end{equation*}
The one-dimensional NURBS basis functions can now be defined by
\begin{equation*}
	R_i(\xi) = \frac{B_{i,\check{p},\vec{t}}(\xi)w_i}{W(\xi)}.
\end{equation*}
The extensions to bivariate NURBS surfaces and trivariate NURBS volumes are straightforward. For NURBS volumes, let $\left\{B_{i_1,\check{p}_1,\vec{t}_1}\right\}_{i_1\in[1,n_1]}$, $\left\{B_{i_2,\check{p}_2,\vec{t}_2}\right\}_{i_2\in[1,n_2]}$ and $\left\{B_{i_3,\check{p}_3,\vec{t}_3}\right\}_{i_3\in[1,n_3]}$ be the sets of B-spline basis functions in $\xi_1$-, $\xi_2$- and $\xi_3$-direction, respectively. These sets have their own order ($\check{p}_1$, $\check{p}_2$ and $\check{p}_3$, respectively) and knot vectors ($\vec{t}_1$, $\vec{t}_2$ and $\vec{t}_3$, respectively). The trivariate NURBS basis functions are then defined by
\begin{equation}\label{Eq:NURBS3D}
	R_{i_1,i_2,i_3}(\vec{\xi}) = \frac{B_{i_1,\check{p}_1,\vec{t}_1}(\xi_1)B_{i_2,\check{p}_2,\vec{t}_2}(\xi_2)B_{i_3,\check{p}_3,\vec{t}_3}(\xi_3)w_{i_1,i_2,i_3}}{W(\vec{\xi})},\qquad \vec{\xi} = [\xi_1,\xi_2,\xi_3]^\transpose
\end{equation}
where the weighting function is now given by
\begin{equation*}
	W(\vec{\xi}) = \sum_{\tilde{i}_1=1}^{n_1}\sum_{\tilde{i}_2=1}^{n_2}\sum_{\tilde{i}_3=1}^{n_3}B_{\tilde{i}_1,\check{p}_1,\vec{t}_1}(\xi_1)B_{\tilde{i}_2,\check{p}_2,\vec{t}_2}(\xi_2)B_{\tilde{i}_3,\check{p}_3,\vec{t}_3}(\xi_3)w_{\tilde{i}_1,\tilde{i}_2,\tilde{i}_3}.
\end{equation*}
A NURBS patch can be represented by the transformation
\begin{equation}
    \vec{X}:[0,1]^3\to\Omega_{\mathrm{b}}\subset\R^3,\quad \vec{\xi} \mapsto \sum_{i_1=1}^{n_1}\sum_{i_2=1}^{n_2}\sum_{i_3=1}^{n_3} R_{i_1,i_2,i_3}(\vec{\xi})\vec{P}_{i_1,i_2,i_3}
\end{equation}
with $\vec{P}_{i_1,i_2,i_3}$ being the control points of the patch. For brevity the PML formulation will be derived for a single NURBS patch representing $\Omega_{\mathrm{b}}$, but the generalization to a multipatch representation is straight forward. Without loss of generalization, we have here assumed normalized knots in the patch in which we want the PML transformation. 


Whenever $\Gamma_{\mathrm{a}}$ is a smooth convex surface we can construct the NURBS patches representing the PML by first finding a NURBS parametrization $\vec{X}_{\mathrm{b}}(\xi_1,\xi_2)$ which is a distance $t_{\textsc{pml}}$ away from $\vec{X}_{\mathrm{a}}(\xi_1,\xi_2)$ (achieved by minimizing $\vec{X}_{\mathrm{b}} - \vec{X}_{\mathrm{a}} - t_{\textsc{pml}}\vec{n}_{\mathrm{a}}$ where $\vec{n}_{\mathrm{a}}$ is the normal vector at $\vec{X}_{\mathrm{a}}$) and then computing a linear lofting between $\vec{X}_{\mathrm{a}}$ and $\vec{X}_{\mathrm{b}}$ to obtain the volumetric NURBS patches.



Starting from a boundary representation of the scatterer $\Gamma$ assume now that a surface-to-volume parametrization, $\vec{X}$, has been found for $\Omega_{\mathrm{a}}$ such that we have a NURBS representation of $\vec{X}_{\mathrm{a}}$ through $\vec{X}_{\mathrm{a}}=\vec{X}\vert_{\Gamma_{\mathrm{a}}}$. That is, we know the control points $\vec{P}_{\mathrm{a},i_1,i_2}$ and can compute the normal vector, $\vec{n}_{\mathrm{a}}$, on $\Gamma_{\mathrm{a}}$ (assuming, without loss of generality, that it is well oriented) through
\begin{equation*}
	\vec{X}_{\mathrm{a}}(\xi_1,\xi_2) = \sum_{i_1=1}^{n_1} \sum_{i_2=1}^{n_2} R_{i_1,i_2}(\xi_1,\xi_2)\vec{P}_{\mathrm{a},i_1,i_2},\qquad \vec{n}_{\mathrm{a}} = \frac{1}{c_0}\pderiv{\vec{X}_{\mathrm{a}}}{\xi_1}\times\pderiv{\vec{X}_{\mathrm{a}}}{\xi_2},\quad c_0 = \left\|\pderiv{\vec{X}_{\mathrm{a}}}{\xi_1}\times\pderiv{\vec{X}_{\mathrm{a}}}{\xi_2}\right\|_2.
\end{equation*}
Find $\vec{X}_{\mathrm{b}}(\xi_1,\xi_2)$ from the following minimization problem
\begin{equation}\label{Eq:minimizationProblem}
	\min_{\vec{P}_{\mathrm{b},i_1,i_2}} \int_{\Gamma_{\mathrm{a}}}\left\|\vec{X}_{\mathrm{b}} - \vec{X}_{\mathrm{a}} - t_{\textsc{pml}}\vec{n}_{\mathrm{a}}\right\|_2^2\idiff \Omega
\end{equation}
where $\vec{n}_{\mathrm{a}}$ is the outward point normal vector at $\Gamma_{\mathrm{a}}$ and
\begin{align*}
	\vec{X}_{\mathrm{b}}(\xi_1,\xi_2) = \sum_{i_1=1}^{n_1} \sum_{i_2=1}^{n_2} R_{i_1,i_2}(\xi_1,\xi_2)\vec{P}_{\mathrm{b},i_1,i_2}.
\end{align*}
To find the minima we differentiate \Cref{Eq:minimizationProblem} w.r.t. the components, $P_{\mathrm{b},i_1,i_2,i}$ (for $i=1,2,3$), of $\vec{P}_{\mathrm{b},i_1,i_2}$ and set this expression to zero
\begin{align}
	&\int_{\Gamma_{\mathrm{a}}}2 R_{i_1,i_2}(\xi_1,\xi_2)\vec{X}_{\mathrm{b}} -2R_{i_1,i_2}(\xi_1,\xi_2)(\vec{X}_{\mathrm{a}} + t_{\textsc{pml}}\vec{n}_{\mathrm{a}})\idiff \Omega = 0,\\
	\Rightarrow &\sum_{j_1=1}^{n_1} \sum_{j_2=1}^{n_2} P_{\mathrm{b},j_1,j_2,j}\int_{\Gamma_{\mathrm{a}}}R_{i_1,i_2}(\xi_1,\xi_2)R_{j_1,j_2}(\xi_1,\xi_2)\idiff \Omega = \int_{\Gamma_{\mathrm{a}}}R_{i_1,i_2}(\xi_1,\xi_2)(\vec{X}_{\mathrm{a}} + t_{\textsc{pml}}\vec{n}_{\mathrm{a}})\idiff \Omega\label{Eq:LinSysEqMinima}
\end{align}
for $i_1=1,\dots,n_1$, $i_2=1,\dots,n_2$, $i=1,2,3$. Which results in a system of equations on the form $\vec{M}\vec{P} = \vec{b}$ where $\vec{M}$ is the mass matrix, $\vec{P}$ contains the coefficients $P_{\mathrm{b},i_1,i_2,i}$ and $\vec{b}$ is the ``force'' vector formed from the right hand side of~\Cref{Eq:LinSysEqMinima}. Note that we could in principle used different weights $w_{i_1,i_2}$ in the NURBS representation of $\vec{X}_{\mathrm{b}}$ relative to that of $\vec{X}_{\mathrm{a}}$. However, this would result in a non-linear problem, and we here choose to use the same weights for both surfaces for simplicity.

As we in the limit $n_{\mathrm{dofs}}\to\infty$ we have $\vec{X}_{\mathrm{b}} = \vec{X}_{\mathrm{a}} + t_{\textsc{pml}}\vec{n}_{\mathrm{a}}$ we can state that this method converges to a representation that for each point on $\Gamma_{\mathrm{b}}$ is a distance $t_{\textsc{pml}}$ normally directed from $\Gamma_{\mathrm{a}}$. The PML layer is then identical to the conformal PML (see~\cite{Beriot2020aap} for details). A simple linear lofting between $\vec{X}_{\mathrm{a}}$ and $\vec{X}_{\mathrm{b}}$ then yields the volumetric representation of the PML layer\footnote{This linear approach of achieving the volumetric PML layer is chosen throughout this work and is arguably the simplest and most rigorous approach.}. That is, the NURBS parametrization of the PML is given by
\begin{equation}\label{Eq:X_lofted}
	\vec{X}(\vec{\xi}) = (1-\xi_3)\vec{X}_{\mathrm{a}}(\xi_1,\xi_2) + \xi_3\vec{X}_{\mathrm{b}}(\xi_1,\xi_2) = \vec{X}_{\mathrm{a}}(\xi_1,\xi_2) + \xi_3(\vec{X}_{\mathrm{b}}(\xi_1,\xi_2)-\vec{X}_{\mathrm{a}}(\xi_1,\xi_2)).
\end{equation}
With some regularity assumptions we can show that this third parametric direction (the absorption direction) is normal not just to $\Gamma_{\mathrm{a}}$ but also $\Gamma_{\mathrm{b}}$. This is important in order to avoid exponential growth of waves of grazing incidence that is improperly aligned with the absorption direction. That is, it is important for waves not to have any directional component directed opposite to the absorption direction. More precisely, if a wave locally has direction vector $\vec{k}$ and the absorption direction at the same location is $\vec{d}_3$, then we must have $\vec{k}\cdot\vec{d}_3 > 0$. This should hold throughout the domain $\Omega_{\mathrm{b}}$. An outline of the argument for why $\vec{n}_{\mathrm{b}}=\vec{n}_{\mathrm{a}}$ goes as follows.

Assume that $\vec{n}_{\mathrm{a}}$ is differentiable w.r.t. $\xi_i$, $i=1,2$. As $\|\vec{n}_{\mathrm{a}}\| = 1$ we have $\pderiv{\|\vec{n}_{\mathrm{a}}\|^2}{\xi_i} = 2\vec{n}_{\mathrm{a}}\cdot \pderiv{\vec{n}_{\mathrm{a}}}{\xi_i}$. That is, both $\pderiv{\vec{n}_{\mathrm{a}}}{\xi_1}$ and $\pderiv{\vec{n}_{\mathrm{a}}}{\xi_2}$ are tangential vectors on $\Gamma_{\mathrm{a}}$. Thus, $\pderiv{\vec{n}_{\mathrm{a}}}{\xi_1}\times\pderiv{\vec{n}_{\mathrm{a}}}{\xi_2} = c_1 \vec{n}_{\mathrm{a}}$ for some function $c_1$. As $\pderiv{\vec{X}_{\mathrm{a}}}{\xi_1}$ and  $\pderiv{\vec{X}_{\mathrm{a}}}{\xi_2}$ are tangent vectors on $\Gamma_{\mathrm{a}}$ as well, $\pderiv{\vec{n}_{\mathrm{a}}}{\xi_1}$ and $\pderiv{\vec{n}_{\mathrm{a}}}{\xi_2}$ can locally be represented as a linear combination of these. That is,
\begin{align*}
	\pderiv{\vec{n}_{\mathrm{a}}}{\xi_1} &= c_{11}\pderiv{\vec{X}_{\mathrm{a}}}{\xi_1} + c_{12}\pderiv{\vec{X}_{\mathrm{a}}}{\xi_2}\\
	\pderiv{\vec{n}_{\mathrm{a}}}{\xi_2} &= c_{21}\pderiv{\vec{X}_{\mathrm{a}}}{\xi_1} + c_{22}\pderiv{\vec{X}_{\mathrm{a}}}{\xi_2}
\end{align*}
for some functions $c_{ij}$, $i,j=1,2$. Thus, we have
\begin{align*}
	\pderiv{\vec{n}_{\mathrm{a}}}{\xi_1}\times\pderiv{\vec{X}_{\mathrm{a}}}{\xi_2} &= c_{11}\underbrace{\pderiv{\vec{X}_{\mathrm{a}}}{\xi_1}\times\pderiv{\vec{X}_{\mathrm{a}}}{\xi_2}}_{=c_0\vec{n}_{\mathrm{a}}} + c_{12}\underbrace{\pderiv{\vec{X}_{\mathrm{a}}}{\xi_2}\times\pderiv{\vec{X}_{\mathrm{a}}}{\xi_2}}_{=0} = c_0 c_{11} \vec{n}_{\mathrm{a}}\\
	\pderiv{\vec{X}_{\mathrm{a}}}{\xi_1}\times\pderiv{\vec{n}_{\mathrm{a}}}{\xi_2} &= c_{21}\underbrace{\pderiv{\vec{X}_{\mathrm{a}}}{\xi_1}\times\pderiv{\vec{X}_{\mathrm{a}}}{\xi_1}}_{=0} + c_{22}\underbrace{\pderiv{\vec{X}_{\mathrm{a}}}{\xi_1}\times\pderiv{\vec{X}_{\mathrm{a}}}{\xi_2}}_{=c_0\vec{n}_{\mathrm{a}}} = c_0 c_{22} \vec{n}_{\mathrm{a}}.
\end{align*}
We can now compute
\begin{align*}
	\pderiv{\vec{X}_{\mathrm{b}}}{\xi_1}\times \pderiv{\vec{X}_{\mathrm{b}}}{\xi_2} &= \left(\pderiv{\vec{X}_{\mathrm{a}}}{\xi_1} + t_{\textsc{pml}}\pderiv{\vec{n}_{\mathrm{a}}}{\xi_1}\right)\times\left(\pderiv{\vec{X}_{\mathrm{a}}}{\xi_2} + t_{\textsc{pml}}\pderiv{\vec{n}_{\mathrm{a}}}{\xi_2}\right)\\
	&= \left[c_0 + t_{\textsc{pml}}c_0(c_{11}+c_{22}) + t_{\textsc{pml}}^2c_1\right]\vec{n}_{\mathrm{a}}
\end{align*}
which implies
\begin{align*}
	\vec{n}_{\mathrm{b}} = \frac{\pderiv{\vec{X}_{\mathrm{b}}}{\xi_1}\times \pderiv{\vec{X}_{\mathrm{b}}}{\xi_2}}{\left\|\pderiv{\vec{X}_{\mathrm{b}}}{\xi_1}\times \pderiv{\vec{X}_{\mathrm{b}}}{\xi_2}\right\|_2} = \vec{n}_{\mathrm{a}}.
\end{align*}
For some boundaries represented by $\vec{X}_{\mathrm{a}}$ the procedure for finding $\vec{X}_{\mathrm{b}}$ satisfy the equality $\vec{X}_{\mathrm{b}} = \vec{X}_{\mathrm{a}} + t_{\textsc{pml}}\vec{n}_{\mathrm{a}}$ exactly, even on the coarsest mesh. This includes spherical, cylindrical, and Cartesian PML layers. The former two are satisfied due to the exact geometry representation of conic sections provided by IGA. For other boundaries we may have (for a finite $n_{\mathrm{dofs}}$) only the approximation $\vec{X}_{\mathrm{b}} \approx \vec{X}_{\mathrm{a}} + t_{\textsc{pml}}\vec{n}_{\mathrm{a}}$. Thus, the decay direction being normal to the boundaries $\Gamma_{\mathrm{a}}$ and $\Gamma_{\mathrm{b}}$ is not exact (only approximated), but this is also the case for the classical FEM formulations provided for PML for non-trivial boundaries. To illustrate this, consider $\Gamma_{\mathrm{a}}$ to be an ellipsoid
\begin{equation}
	\Gamma_{\mathrm{a}} = \left\{\vec{x}=\begin{bmatrix}
		x_1\\x_2\\x_3
	\end{bmatrix}\in\R^3\st \left(\frac{x_1}{a_1}\right)^2 + \left(\frac{x_2}{a_2}\right)^2 + \left(\frac{x_3}{a_3}\right)^2 = 1\right\}
\end{equation}
where $a_i$ is the semi major/minor axis of the ellipsoid in the Cartesian direction $i$. The spherical parametrization can be extended for this with the parametrization
\begin{equation}\label{Eq:X_a_trig}
	\vec{X}_{\mathrm{a}}(\vartheta,\varphi) = \begin{bmatrix}
		a_1 \sin\vartheta\cos\varphi\\
		a_2 \sin\vartheta\sin\varphi\\
		a_3 \cos\vartheta
	\end{bmatrix},\quad \vartheta\in[0,\pi]\quad\varphi\in[0,2\pi].
\end{equation}
From this we find
\begin{equation}\label{Eq:X_b_trig}
	\vec{X}_{\mathrm{b}}(\vartheta,\varphi) = \vec{X}_{\mathrm{a}}(\vartheta,\varphi) +  \frac{t}{q(\vartheta,\varphi)}\begin{bmatrix}
		a_2a_3\sin\vartheta\cos\varphi\\
		a_1a_3\sin\vartheta\sin\varphi\\
		a_1a_2\cos\vartheta
	\end{bmatrix}
\end{equation}
where
\begin{equation*}
	q(\vartheta,\varphi) = \sqrt{a_1^2 a_2^2 \cos^2\vartheta + \sin^2\vartheta\left(a_1^2 a_3^2 \sin^2\varphi + a_2^2 a_3^2 \cos^2\varphi\right)}.
\end{equation*}
NURBS representations of $\vec{X}_{\mathrm{a}}$ exists (c.f.~\cite{Venas2015iao}) but it is not clear if it exists for $\vec{X}_{\mathrm{b}}$ (at least it requires some work to find) and so we need to approximate this surface using IGA.

\begin{figure}
	\centering
	\begin{subfigure}{0.41\textwidth}
		\centering
		\includegraphics[width=\textwidth]{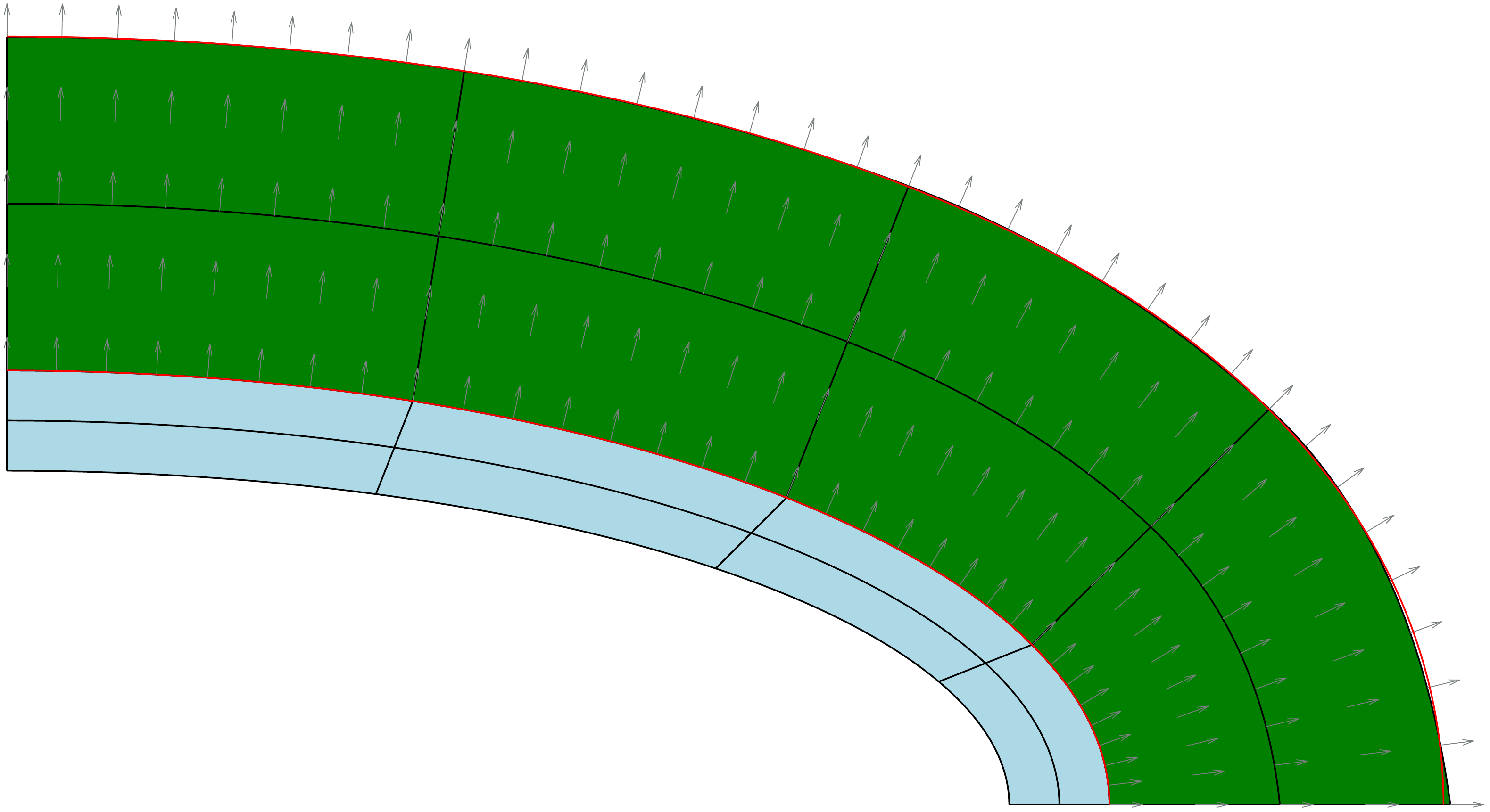}
		\caption{An IGA mesh with $p=2$ and maximum continuity and $\Gamma_{\mathrm{b}}$ obtained by least squares approach. The number of degrees of freedom is 1092.}
		\label{Fig:PML_IGA}
    \end{subfigure}
    \hspace{1cm}
	\begin{subfigure}{0.41\textwidth}
		\centering
		\includegraphics[width=\textwidth]{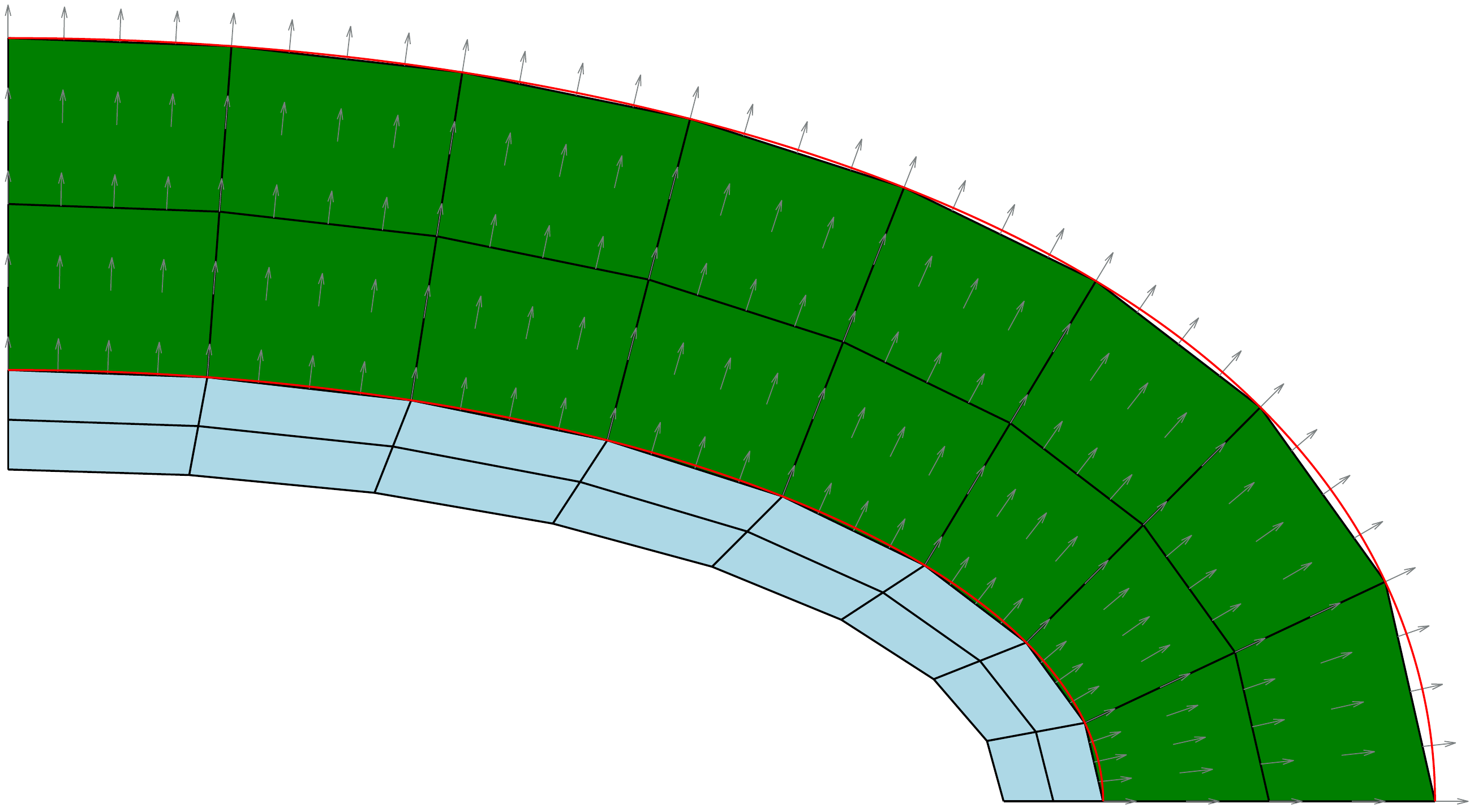}
		\caption{A classical FEM linear discretization and $\Gamma_{\mathrm{b}}$ obtained by interpolation. The number of degrees of freedom is 1928.}
		\label{Fig:PML_linear_FEM}
    \end{subfigure}
    \caption{A cross section of the meshes around a ellipsoidal scatterer in the first quadrant of the $xy$-plane. The red curves are the exact trigonometric parametrization of $\vec{X}_{\mathrm{a}}$ and $\vec{X}_{\mathrm{b}}$ in~\Cref{Eq:X_a_trig,Eq:X_b_trig}, respectively. The grey arrows show the direction of the absorption in the (green) PML layer.}
	\label{Fig:PMLillustration}
\end{figure}

\begin{figure}
	\centering
	\begin{subfigure}{0.4\textwidth}
		\centering
		\includegraphics[width=\textwidth]{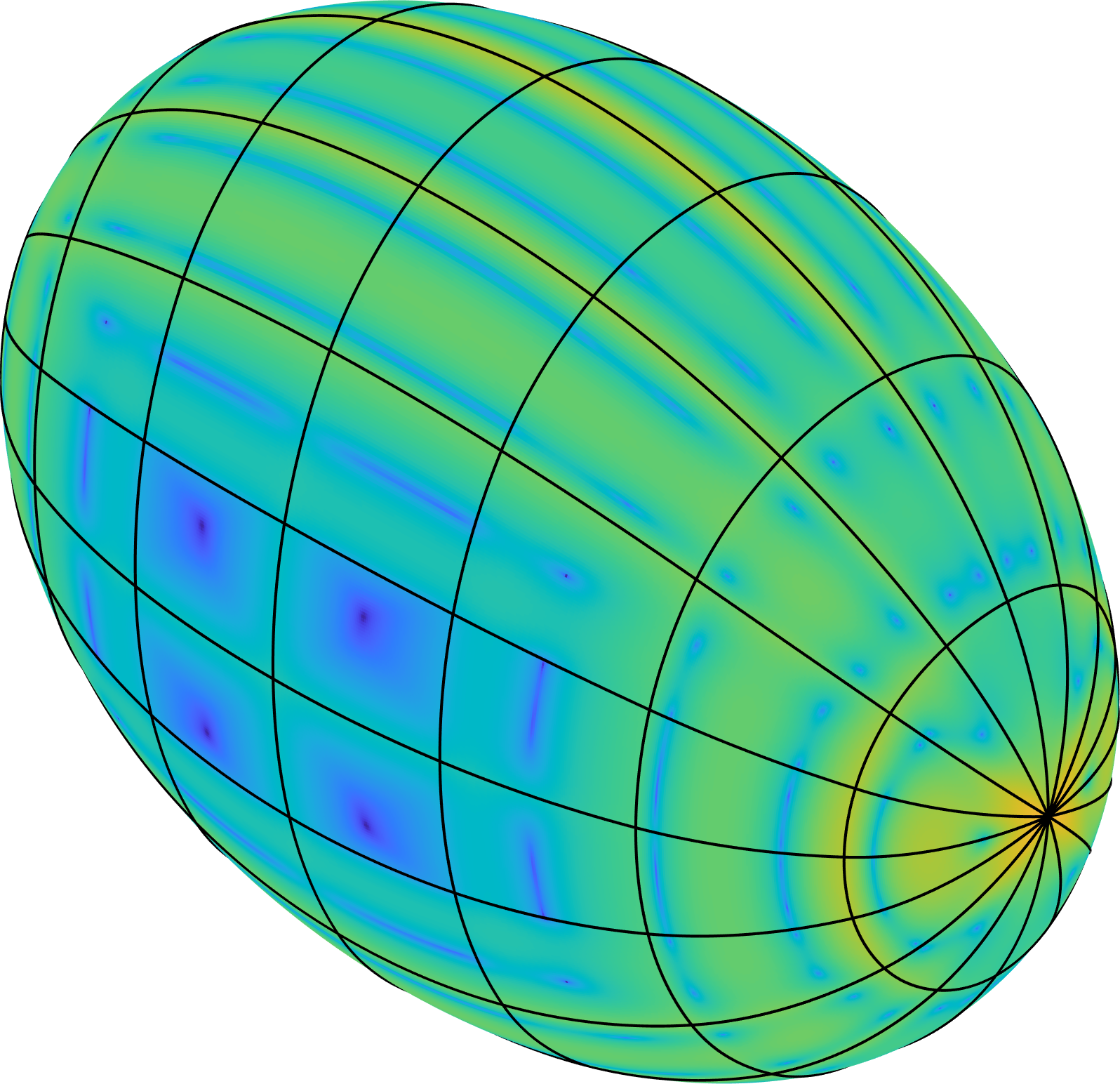}
		\caption{An IGA mesh with $\check{p}=2$ and maximum continuity for $\Gamma_{\mathrm{b}}$ obtained by least squares approach. The number of degrees of freedom is 1092.}
		\label{Fig:PML_IGA_3D}
    \end{subfigure}
	\begin{subfigure}{0.15\textwidth}
		\centering
		\includegraphics[width=\textwidth]{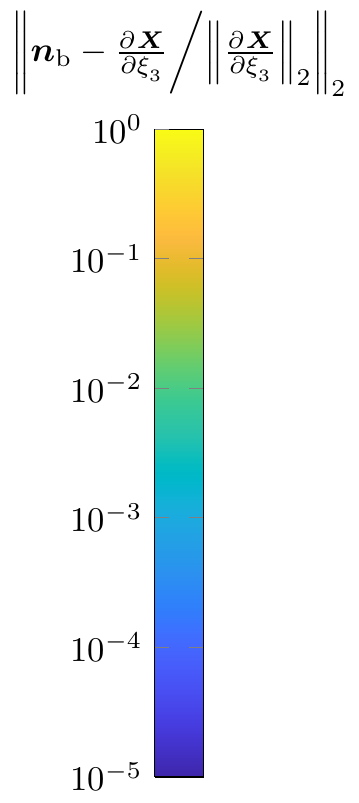}
    \end{subfigure}
	\begin{subfigure}{0.4\textwidth}
		\centering
		\includegraphics[width=\textwidth]{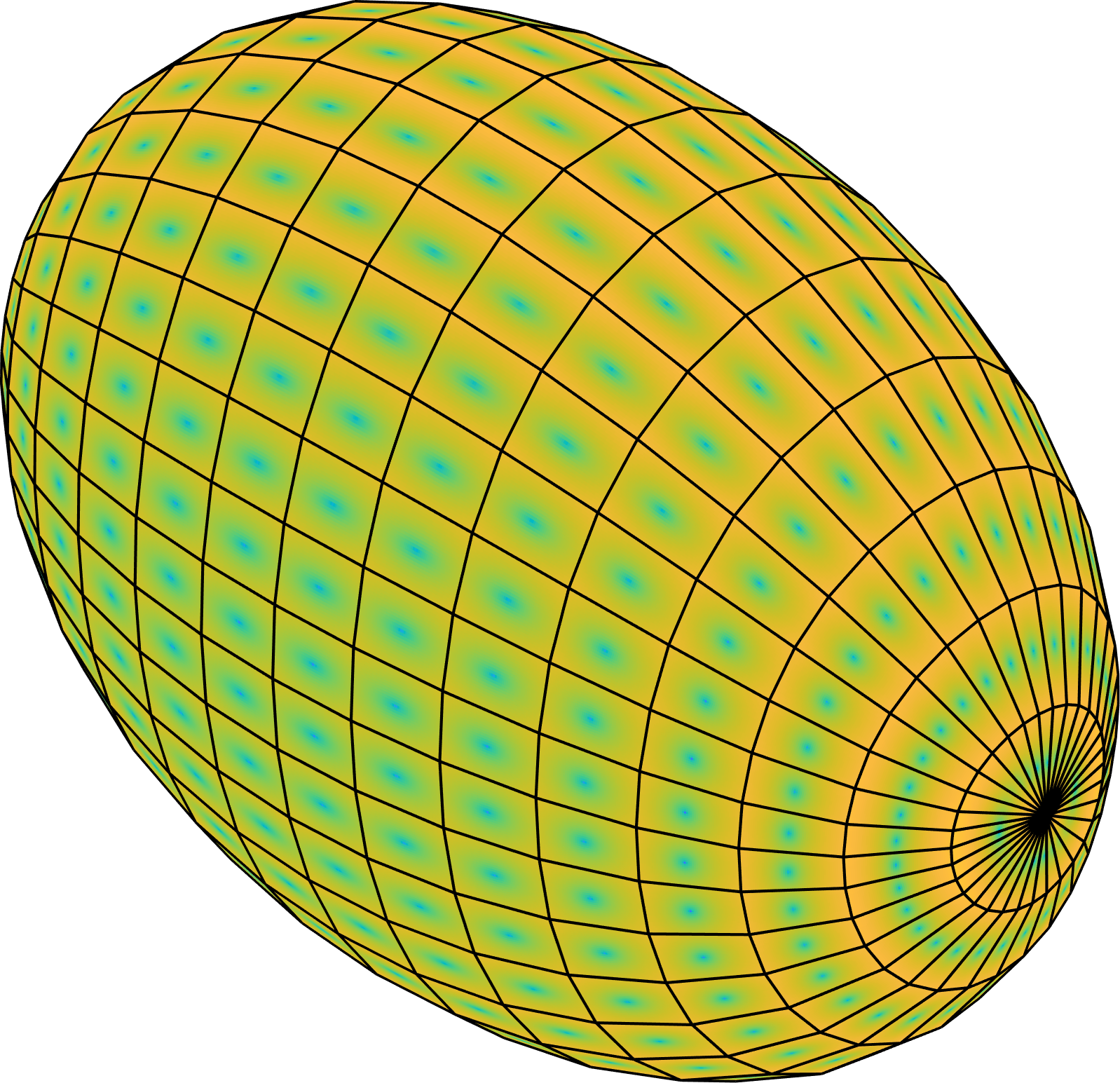}
		\caption{A classical FEM linear discretization for $\Gamma_{\mathrm{b}}$ obtained by interpolation. The number of degrees of freedom is 1928.}
		\label{Fig:PML_linear_FEM_3D}
    \end{subfigure}
    \caption{The boundary of $\Gamma_{\mathrm{b}}$ approximated with IGA and linear FEM. The coloring is made with the error between the normal vector and the absorption direction. In the $L^2$-norm this error computes to 1.03\% and 7.16\% for IGA and linear FEM, respectively, even though linear FEM here uses 1.77 times as many dofs.}
	\label{Fig:PMLillustration_3D}
\end{figure}

\begin{figure}
	\centering
	\includegraphics[width=\textwidth]{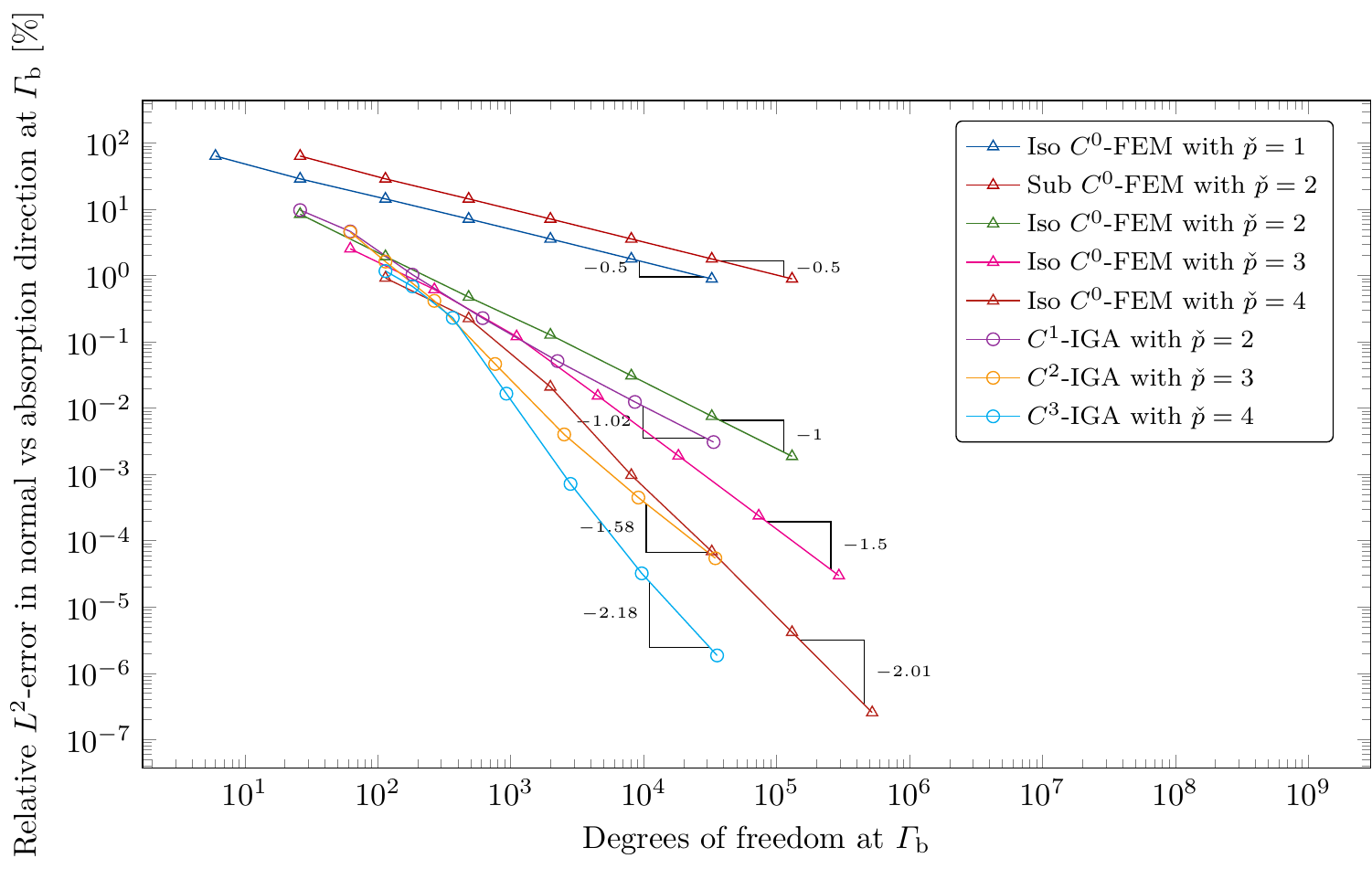}
    \caption{The relative $L^2$-error between the normal vector, $\vec{n}$, and the absorption direction $\pderiv{\vec{X}}{\xi_3}\bigg/\left\|\pderiv{\vec{X}}{\xi_3}\right\|_2$ is plotted against the number of degrees of freedom at $\Gamma_{\mathrm{b}}$. Notice the ``left shift'' obtained for the $C^{\check{p}-1}$-IGA-curves compared to the $C^0$-FEM-curves due to higher regularity which implies that $C^{\check{p}-1}$-IGA is more accurate per degrees of freedom than $C^0$-FEM. This effect increases with increasing polynomial order. Furthermore, the use of {\it Iso}-parametric FEM is as expected superior to {\it Sub}-parametric FEM.}
	\label{Fig:PML_normalConvergence}
\end{figure}

Consider an ellipsoidal scatterer with $a_1=\SI{3}{m}$, $a_2=\SI{1}{m}$ and $a_3=\SI{2}{m}$. We construct the artificial boundary $\Gamma_{\mathrm{a}}$ through a larger ellipsoidal with $a_1=\SI{3.3}{m}$, $a_2=\SI{1.3}{m}$ and $a_3=\SI{2.3}{m}$. Both of these surfaces can be exactly parametrized with NURBS and so we find the volumetric parametrization for $\Omega_{\mathrm{a}}$ with a simple linear lofting between $\Gamma$ and $\Gamma_{\mathrm{a}}$. Note that this does not yields ``radial'' mesh lines normal to the surfaces. We now investigate the construction of a PML layer with thickness $t_{\textsc{pml}}$. As we can see from \Cref{Fig:PMLillustration}, IGA is able to exactly represent $\Gamma$ and $\Gamma_{\mathrm{a}}$ but not $\Gamma_{\mathrm{b}}$. However, the latter is still much better approximated using IGA compared to linear FEM even with much less degrees of freedom used. This results in a better approximation of mesh lines being normal to the boundaries compared to that of linear FEM.
This is illustrated even better in~\Cref{Fig:PML_linear_FEM_3D} where the expected pattern for the linear FEM emerges; the normal vectors at the approximated boundary $\Gamma_{\mathrm{b}}$ will here be equal to the absorption direction $\pderiv{\vec{X}}{\xi_3}$ only in the center of each element. The results for IGA using interpolation at Greville abscissae gave errors roughly in the same range as that of the least squares.
The comparison is illustrated more rigorously in~\Cref{Fig:PML_normalConvergence}.

The complex coordinate transformation representing the PML coordinate stretching is given by
\begin{equation}\label{Eq:S}
    \vec{S}:[0,1]^3\to \C^3,\quad \vec{\xi} \mapsto \vec{S}(\vec{\xi}) = \tilde{\vec{\xi}},\qquad \tilde{\vec{\xi}} = [\tilde{\xi}_1,\tilde{\xi}_2,\tilde{\xi}_3]
\end{equation}
where
\begin{equation}
	\tilde{\xi}_i = \xi_i+\imag I_i(\xi_i),\quad I_i(\xi_i) = \int_0^{\xi_i}\sigma_i(\xi)\idiff\xi
\end{equation}
and $\sigma_i$ is a monotonically increasing function satisfying $\sigma_i\geq 0$. 

Instead of using the bilinear form defined in the physical space ($\vec{x}\in\Omega^+\subset\R^3$), we use the bilinear form over the space $\tilde{\Omega}^+ = \{\vec{T}(\vec{x})\st \vec{x}\in\Omega^+\}$ where
\begin{equation}
    \vec{T}: \R^3 \to \C^3,\quad \tilde{\vec{x}} = \vec{T}(\vec{x}) = \vec{X}(\vec{S}(\vec{X}^{-1}(\vec{x}))).
\end{equation}
The bilinear form is then given by
\begin{equation}
    B(p,q) = \int_{\tilde{\Omega}^+}\tilde{\nabla}q\cdot\tilde{\nabla}p - k^2qp\idiff\tilde{\Omega}.
\end{equation}
We then need the Jacobian matrix, $\tilde{\vec{J}}=\pderiv{\tilde{\vec{x}}}{\vec{x}}$, in order to compute $\tilde{\nabla} p$ from $\tilde{\nabla} p\tilde{\vec{J}} = \nabla p$. Applying the chain rule (with $\vec{J} =\pderiv{\vec{X}}{\vec{\xi}}$ and $\vec{D}=\pderiv{\tilde{\vec{\xi}}}{\vec{\xi}} = \vec{I}+\imag\diag{\sigma_1(\xi_1), \sigma_2(\xi_2), \sigma_3(\xi_3)}$) yields
\begin{equation}\label{Eq:J_tilde}
    \tilde{\vec{J}} = \pderiv{\tilde{\vec{X}}}{\tilde{\vec{\xi}}}\pderiv{\tilde{\vec{\xi}}}{\vec{\xi}}\pderiv{\vec{\xi}}{\vec{X}} = \vec{J}(\tilde{\vec{\xi}})\vec{D}(\vec{\xi})\vec{J}^{-1}(\vec{x}).
\end{equation}
Since $\nabla p \vec{J}(\vec{\xi}) = \nabla_{\vec{\xi}} p$ we have $\tilde{\nabla} p = \nabla_{\vec{\xi}} p\vec{D}^{-1}\vec{J}(\tilde{\vec{\xi}})^{-1}$, which inserted into the bilinear form yields (for a single patch)
\begin{equation}
    B(p,q) = \int_{[0,1]^3}\left[\left(\vec{J}(\tilde{\vec{\xi}})^{-\transpose}\vec{D}^{-1}\nabla_{\vec{\xi}}^\transpose q\right)\cdot\left(\vec{J}(\tilde{\vec{\xi}})^{-\transpose}\vec{D}^{-1}\nabla_{\vec{\xi}}^\transpose p\right) - k^2qp\right]\det\left(\vec{J}(\tilde{\vec{\xi}})\right)\det\left(\vec{D}\right)\idiff\vec{\xi}.
\end{equation}
This reduces to the standard bilinear form whenever $\sigma_i=0,\,\forall i$. Compare this bilinear form to the bilinear form for spherical coordinates, $(r,\vartheta,\varphi)$ in~\cite{Shirron2006afe}
\begin{equation}
    B(p,q) = \int_{[0,1]^3}\left[(\vec{J}_{\mathrm{s}}^{-\transpose}\vec{D}^{-1}\vec{J}_{\mathrm{s}}^\transpose\vec{J}^{-\transpose}\nabla_{\vec{\xi}}^\transpose q)\cdot(\vec{J}_{\mathrm{s}}^{-\transpose}\vec{D}^{-1}\vec{J}_{\mathrm{s}}^\transpose\vec{J}^{-\transpose}\nabla_{\vec{\xi}}^\transpose p) - k^2qp\right]\det(\vec{J})\det(\vec{D})\idiff\vec{\xi}
\end{equation}
where
\begin{equation}
	\vec{J}_{\mathrm{s}} = \begin{bmatrix}
		\sin\vartheta\cos\varphi & r\cos\vartheta\cos\varphi & -r\sin\vartheta\sin\varphi\\
		\sin\vartheta\sin\varphi & r\cos\vartheta\sin\varphi & r\sin\vartheta\cos\varphi\\
		\cos\vartheta & -r\sin\vartheta & 0
	\end{bmatrix},\quad
    \vec{D} = \vec{I}+\imag\diag{\sigma(\xi), \frac{1}{r}I(\xi), \frac{1}{r}I(\xi)}
\end{equation}
with $\xi=\frac{r-R}{S-R}$. In the present formulation the evaluation of trigonometric functions in the assembly procedure is therefore replaced by NURBS-evaluations.

However, the formula for the Jacobian in~\Cref{Eq:J_tilde} requires NURBS evaluations with complex parametric argument. For linear absorption parametrizations (e.g.~\Cref{Eq:X_lofted}) this can be avoided. Starting by differentiation of~\Cref{Eq:X_lofted} w.r.t. $\xi_3$ yields
\begin{equation}
	\pderiv{\vec{X}}{\xi_3} = \vec{X}_{\mathrm{b}}(\xi_1,\xi_2)-\vec{X}_{\mathrm{a}}(\xi_1,\xi_2)
\end{equation}
which is constant in $\xi_3$ due to the linearity in this parametric direction. 
Assume we now only have a single absorption direction in the third parametric direction. From~\Cref{Eq:X_lofted} and~\Cref{Eq:S} we then have
\begin{align*}
 	\tilde{\vec{x}} &= \vec{X}(\xi_1,\xi_2,\xi_3+\imag I_3(\xi_3)) = \vec{X}(\vec{\xi}) + \imag I_3(\xi_3)(\vec{X}_{\mathrm{b}}(\xi_1,\xi_2)-\vec{X}_{\mathrm{a}}(\xi_1,\xi_2))\\
 	&= \vec{X}(\xi_1,\xi_2,\xi_3) + \imag I_3(\xi_3)\pderiv{\vec{X}}{\xi_3}
\end{align*}
such that
\begin{align*}
 	\pderiv{\tilde{\vec{x}}}{\vec{\xi}} &= \vec{J} + \imag \left[I_3(\xi_3)\ppderiv{\vec{X}}{\xi_1}{\xi_3},I_3(\xi_3)\ppderiv{\vec{X}}{\xi_2}{\xi_3}, I_3'(\xi_3)\pderiv{\vec{X}}{\xi_3}\right]\\
&=   \vec{J} + \imag \left[\ppderiv{\vec{X}}{\xi_1}{\xi_3},\ppderiv{\vec{X}}{\xi_2}{\xi_3}, \pderiv{\vec{X}}{\xi_3}\right]\diag{I_3(\xi_3),I_3(\xi_3),\sigma_3(\xi_3)}.
\end{align*}
If we have absorption in the $i$-th parametric direction, we can with the partition of unity property of NURBS write
\begin{equation}
	\ppderiv{\vec{X}}{\xi_j}{\xi_i} = \sum_{i_1=1}^{n_1}\sum_{i_2=1}^{n_2}\sum_{i_3=1}^{n_3} \frac{1}{R_{i_1,i_2,i_3}(\vec{\xi})}\pderiv{R_{i_1,i_2,i_3}}{\xi_j}\pderiv{R_{i_1,i_2,i_3}}{\xi_i}\vec{P}_{i_1,i_2,i_3},\qquad j\neq i
\end{equation}
for any $\vec{\xi}$ inside an element. 
As the first order derivatives of the NURBS basis functions are readily available in all IGA codes no extra basis function evaluations are needed. The Jacobian expressions for multiple absorption direction are slightly more cumbersome and is outlined in~\Cref{Sec:multiplePMLdirections}. However, the main approach presented herein only have a single absorption direction (in the $\xi_3$-direction) which liberates us from these expressions. The main point here is that there is no need to evaluate NURBS functions with complex parametric arguments (as~\Cref{Eq:S} might indicate). However, it is convenient to have such an extension implemented for code verification. For industrial codes the extension in~\Cref{Eq:orderOneBspline} for these cases is thus not really needed for the linear lofting approach suggested in the present work.
\subsection{Absorption functions}
In~\cite{Shirron2006afe}, the following decay function (or absorption function) is used
\begin{equation}\label{Eq:sigmaType1}
    \sigma(\xi) = \xi(\euler^{\gamma\xi} - 1)
\end{equation}
which gives
\begin{equation}
	I(\xi) = \frac{\euler^{\gamma\xi}(\gamma\xi-1)+1}{\gamma^2}-\frac{\xi^2}{2}.
\end{equation}
This would require finding $\gamma$ for each setup as it would be depending on the frequency and the PML thickness. Alternatively, the following decay function may be used~\cite{Michler2007itp,Astaneh2018opm,Mi2021ilc}
\begin{equation}\label{Eq:sigmaType2}
    \sigma(\xi) = -\gamma\xi^n\ln\epsilon,\quad n=2
\end{equation}
which gives
\begin{equation}
	I(\xi) = -\frac{\gamma}{n+1}\xi^{n+1}\ln\epsilon.
\end{equation}
In~\cite{Bermudez2007aop,Bermudez2008aeb} a decay function with unbounded integral was shown to be optimal with the assumption of Dirichlet boundary conditions at $\Gamma_{\mathrm{b}}$. The following function
\begin{equation}\label{Eq:sigmaType3}
    \sigma(\xi) = \gamma(1-\xi)^{-n},\quad 1\leq n < 3
\end{equation}
with
\begin{equation}
	I(\xi) = \begin{cases}
		\gamma\frac{(1-\xi)^{1-n} -1}{n-1} & n > 1\\
		-\gamma\ln(1-\xi) & n = 1
	\end{cases}
\end{equation}
was found to be optimal for $n=1$ and $\gamma=\frac{1}{k t_{\textsc{pml}}}$ (translated for the present PML formulation) for 2D acoustic scattering. Unless otherwise stated, this function will be used in the examples herein. As noted in~\cite{Bermudez2007aop} this function yields discontinuity at $\xi=0$. Somewhat surprisingly, the continuous alternative
\begin{equation}\label{Eq:sigmaType4}
    \sigma(\xi) = \gamma\left[(1-\xi)^{-n} - 1\right],\quad 1\leq n < 3
\end{equation}
with
\begin{equation}
	I(\xi) = \begin{cases}
		\gamma\left[\frac{(1-\xi)^{1-n} -1}{n-1}-\xi\right] & n > 1\\
		-\gamma\left[\ln(1-\xi)+\xi\right] & n = 1
	\end{cases}
\end{equation}
did not give better results.

Consider a plane wave at the far side of the PML where the PML (between $0<x_3<t_{\textsc{pml}}$) is parameterized by $\vec{X} = \xi_1\vec{e}_1+\xi_2\vec{e}_2+t_{\textsc{pml}}\xi_3\vec{e}_3$ (where $\vec{e}_i$, $i=1,2,3$ are the standard Cartesian basis vectors in $\R^3$)
\begin{equation}
    p_{\mathrm{inc}}(\tilde{\vec{x}})\vert_{\tilde{\Gamma_{\mathrm{b}}}} = P_{\mathrm{inc}}\euler^{\imag \vec{k}\cdot\tilde{\vec{x}}} = P_{\mathrm{inc}}\euler^{\imag \vec{k}\cdot\vec{x}}\euler^{-k_3 t_{\textsc{pml}}I(1)}.
\end{equation}
If we want $|p_{\mathrm{inc}}/P_{\mathrm{inc}}|$ to decay to a value $\epsilon$ at $\Gamma_{\mathrm{b}}$ we must have
\begin{equation}
    \epsilon = \euler^{-k_3 t_{\textsc{pml}}I(1)}.
\end{equation}
For the particular function in \Cref{Eq:sigmaType2} we can compute $\gamma$ to be
\begin{equation}
    \gamma = \frac{n+1}{k_3 t_{\textsc{pml}}}.
\end{equation}
\section{Numerical examples} 
\label{sec:resultsDisc}
We initiate this section with an investigation on a sphere where analytic solution exists to the plane wave scattering problem. Then, we consider a manufactured solution on a cylindrical domain, before ending with a scattering problem on a more complex geometry.

\subsection{Scattering from rigid sphere}
\begin{figure}
	\centering
	\begin{subfigure}{0.3\textwidth}
		\centering
		\includegraphics[width=0.95\textwidth]{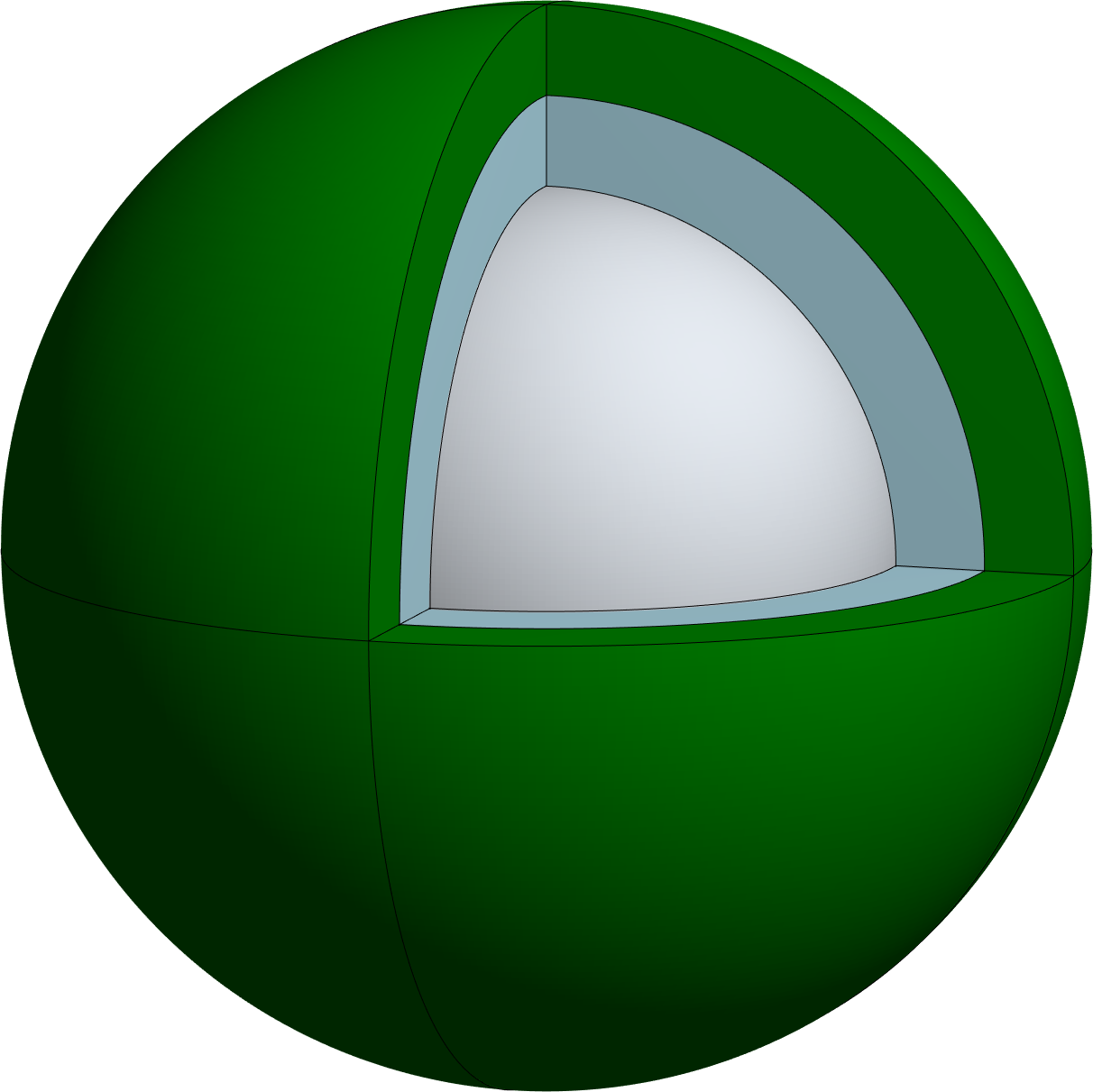}
		\caption{Mesh ${\cal M}_{1,\check{p},\check{k}}^{\textsc{igapml}}$}
		\label{Fig2:SphericalShellMeshes1}
    \end{subfigure}
    ~
	\begin{subfigure}{0.3\textwidth}
		\centering
		\includegraphics[width=0.95\textwidth]{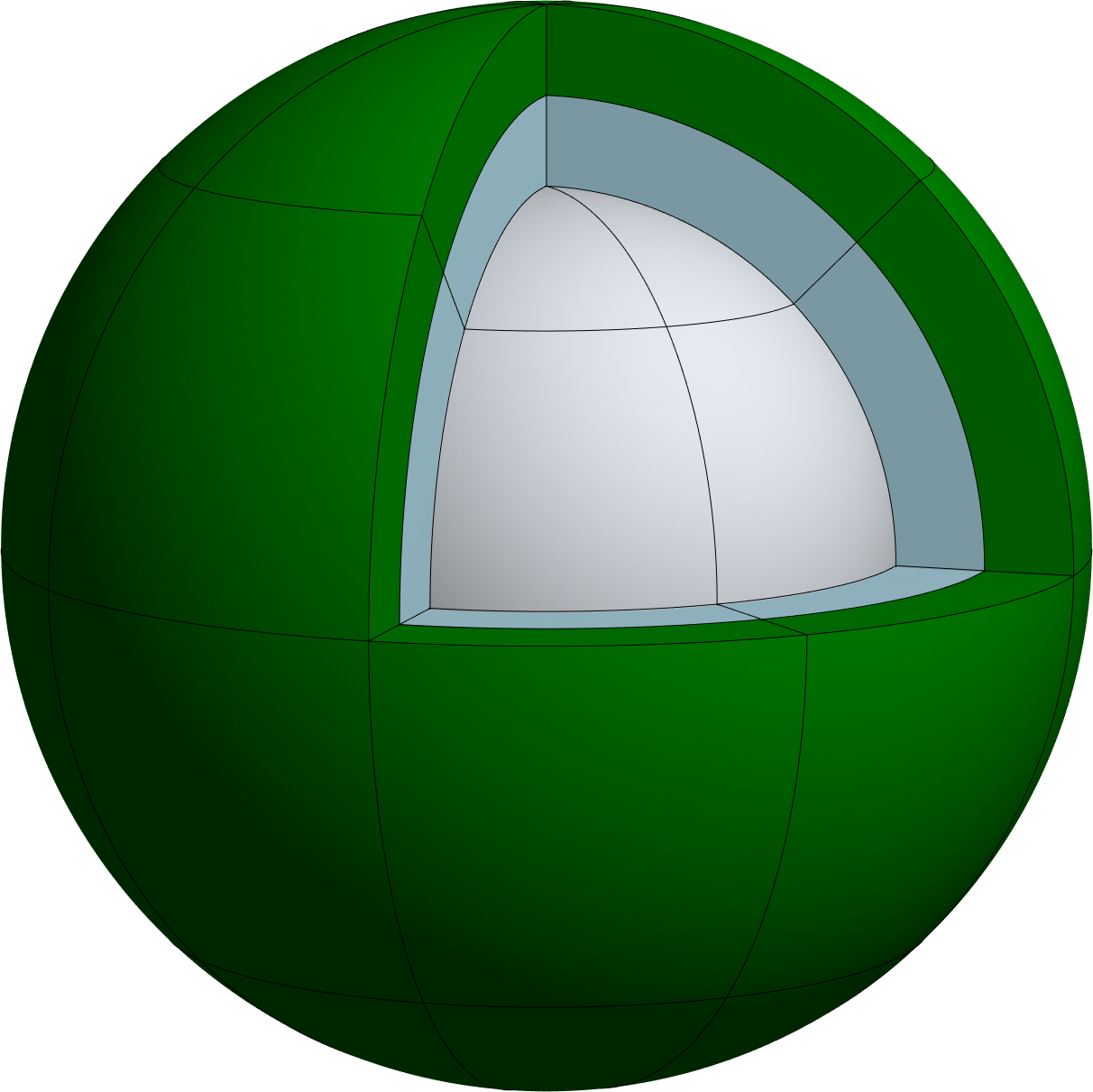}
		\caption{Mesh ${\cal M}_{2,\check{p},\check{k}}^{\textsc{igapml}}$}
		\label{Fig2:SphericalShellMeshes2}
    \end{subfigure}
    ~
	\begin{subfigure}{0.3\textwidth}
		\centering
		\includegraphics[width=0.95\textwidth]{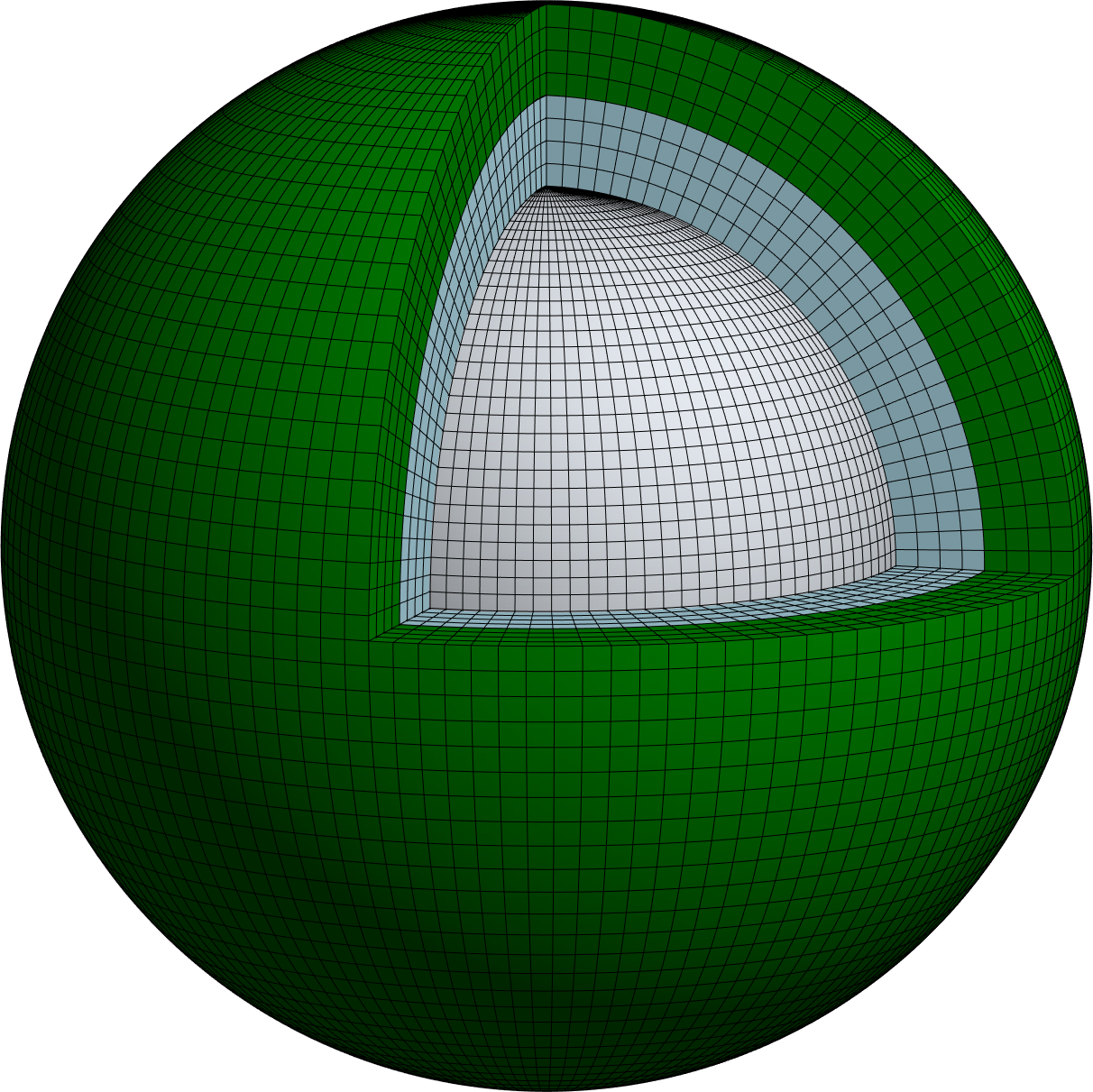}
		\caption{Mesh ${\cal M}_{6,\check{p},\check{k}}^{\textsc{igapml}}$}
		\label{Fig2:SphericalShellMeshes3}
    \end{subfigure}
	\caption{\textbf{Scattering from rigid sphere}: Illustration of the coarse mesh ${\cal M}_{1,\check{p},\check{k}}^{\textsc{igapml}}$ and the first and fifth refinement. The PML domain (green) has the same thickness as the domain inside $\Gamma_{\mathrm{a}}$ (in light blue) which is attached to the spherical (grey) scatterer, $\Gamma$.}
	\label{Fig2:SphericalShellMeshes}
\end{figure}

We start by performing the same analysis done in Figure 9 in~\cite{Venas2018iao} were the convergence through $h$-refinement is studied on a rigid scattering problem on a sphere of radius $R=\SI{5.075}{m}$. For completeness the mesh construction is here repeated. The meshes will be generated from a standard discretization of a sphere using NURBS as seen in \Cref{Fig2:SphericalShellMeshes}. We shall denote by ${\cal M}_{m,\check{p},\check{k}}^{\textsc{igapml}}$, mesh number $m$ with polynomial order $\check{p}$ and continuity $\check{k}$ across element boundaries\footnote{Except for some possible $C^0$ interfaces in the initial CAD geometry.}. For the corresponding FEM meshes we denote by ${\cal M}_{m,\check{p},\mathrm{s}}^{\textsc{fempml}}$ and ${\cal M}_{m,\check{p},\mathrm{i}}^{\textsc{fempml}}$ the subparametric and isoparametric FEM meshes, respectively. The initial mesh is depicted as mesh ${\cal M}_{1,\check{p},\check{k}}^{\textsc{igapml}}$ in \Cref{Fig2:SphericalShellMeshes1} and is refined only in the angular directions for the first~3 refinements (that is, mesh ${\cal M}_{4,\check{p},\check{k}}^{\textsc{igapml}}$ only have two element thickness in the radial direction). Mesh ${\cal M}_{m,\check{p},\check{k}}^{\textsc{igapml}}$, $m=5,6$, have 4 and 8 elements in the radial direction, respectively. This is done to obtain low aspect ratios for the elements. All the meshes will then be nested and the refinements are done uniformly. We shall use the same polynomial order in all parameter directions; $\check{p}_1=\check{p}_2=\check{p}_3$. The ${\cal M}_{m,\check{p},\check{k}}^{\textsc{igaiem}}$ meshes used in~\cite{Venas2018iao} correspond to the light blue domain in~\Cref{Fig2:SphericalShellMeshes} (with the PML-layer being replaced by infinite elements).

The energy norm is here defined by (cf.~\cite{Venas2018iao})
\begin{equation}
    \label{Eq:energyNorm}
	\energyNorm{p}{\Omega_{\mathrm{a}}} = \sqrt{\int_{\Omega_{\mathrm{a}}} \left|\nabla p\right|^2 + k^2|p|^2 \idiff\Omega}.
\end{equation}

\begin{figure}
	\centering
	\includegraphics[width=\textwidth]{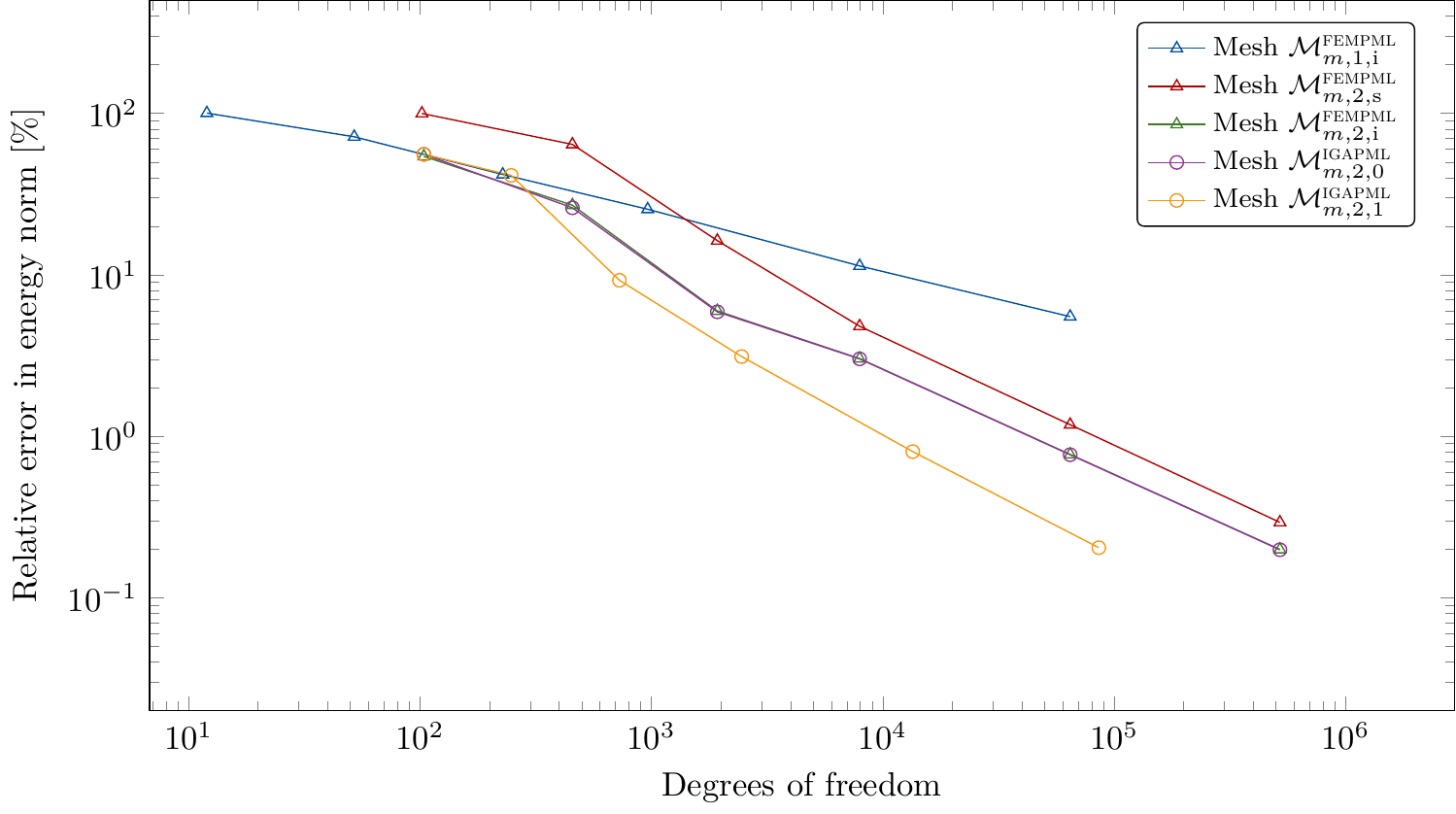}
	\caption{\textbf{Scattering from rigid sphere}: Convergence analysis on the rigid scattering case with $kR=5.075$ and mesh ${\cal M}_{m,2,1}$, $m=1,\dots,6$. For the PML formulation the stretching function in~\Cref{Eq:sigmaType3} is used with $n=1$. The relative energy error (from \Cref{Eq:energyNorm}) is plotted against the degrees of freedom.}
	\label{Fig2:EnergyErrorPlotsDofs}
\end{figure}

\begin{figure}
	\centering
	\includegraphics[width=\textwidth]{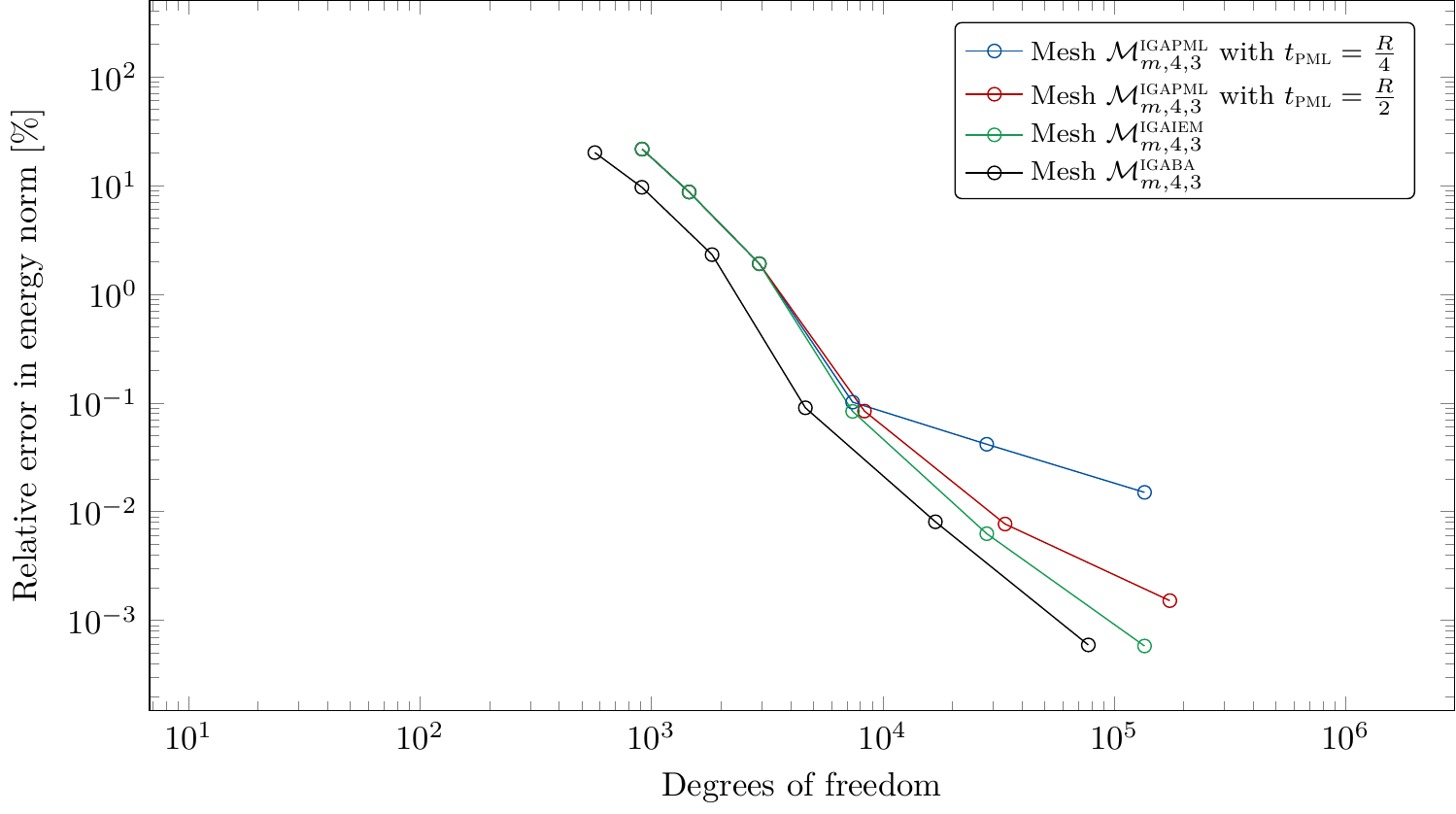}
	\caption{\textbf{Scattering from rigid sphere}: Convergence analysis on the rigid scattering case with $kR=5.075$ and mesh ${\cal M}_{m,4,3}$, $m=1,\dots,6$. For the PML formulation the stretching function in~\Cref{Eq:sigmaType3} is used with $n=1$. The relative energy error (from \Cref{Eq:energyNorm}) is plotted against the degrees of freedom. When doubling the PML thickness we also double the number of elements in the radial direction (which can be seen by an increase of dofs used in the final three meshes). The best approximation (BA) is in the $L^2$-norm over the domain $\Omega_{\mathrm{a}}$. Note that the best approximation simulation lacks elements in $\Omega_{\mathrm{b}}$ and its curve is here thus shifted to the left. With increased PML thickness the PML-simulations converge to the accuracy of the IEM-simulations. }
	\label{Fig2:EnergyErrorPlotsDofs_2}
\end{figure}

\begin{figure}
	\centering
	\includegraphics[width=\textwidth]{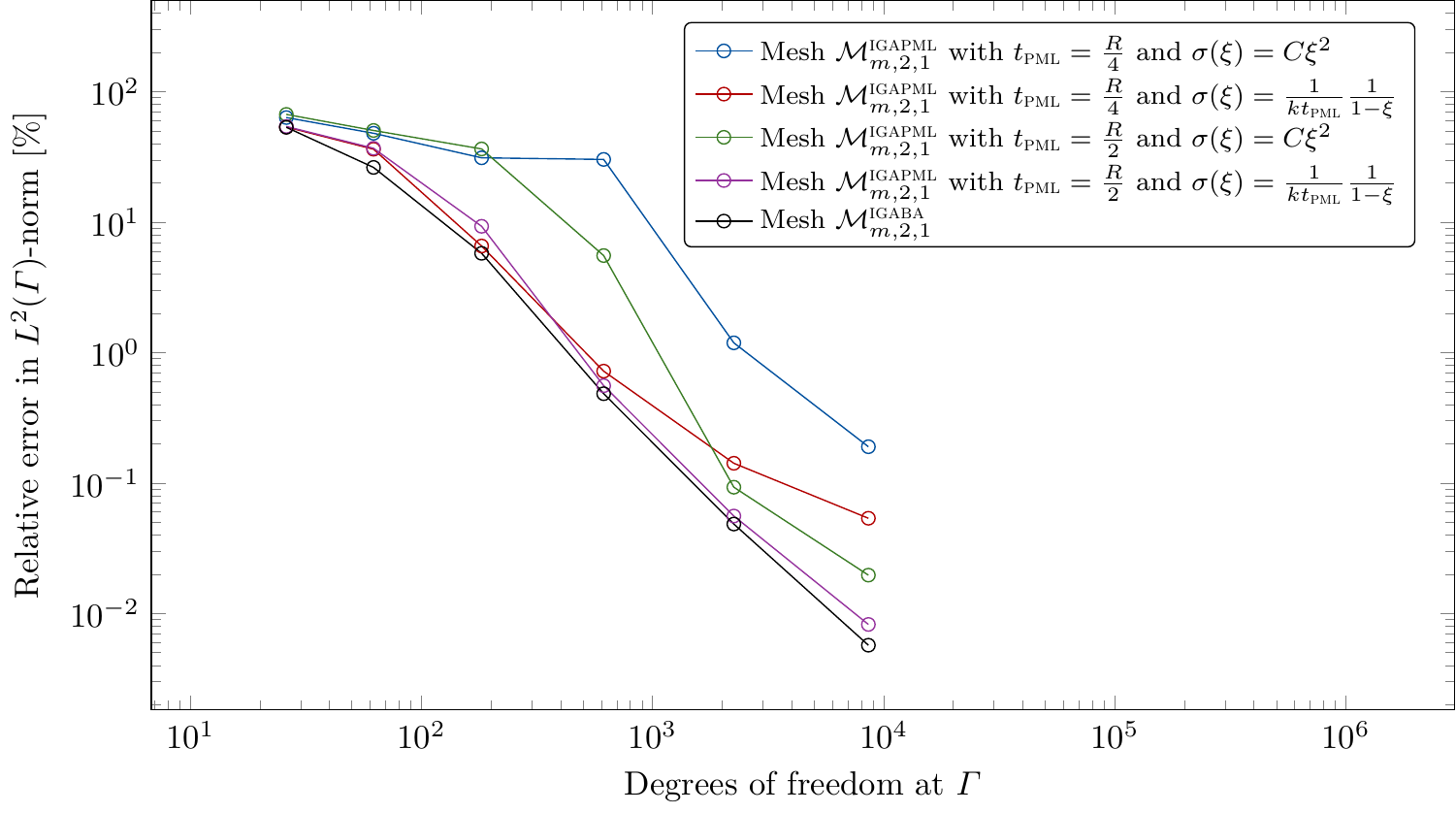}
	\caption{\textbf{Scattering from rigid sphere}: Convergence analysis on the rigid scattering case with $kR=5.075$ and mesh ${\cal M}_{m,4,3}$, $m=1,\dots,6$. The relative $L^2$-error is plotted against the degrees of freedom. The two stretching functions used here are given in~\Cref{Eq:sigmaType2,Eq:sigmaType3}, respectively. Here, $C=\alpha/(knt_{\textsc{pml}}^{n-1})$ with $\alpha=30$ and $n=2$. Note that $\alpha=30$ was used instead of $\alpha=10$ as in~\cite{Mi2021ilc} as the former gave much better results. The best approximation (BA) is for the degrees of freedom at $\Gamma$ in the $L^2$-norm.}
	\label{Fig2:EnergyErrorPlotsDofs_3}
\end{figure}

\Cref{Fig2:EnergyErrorPlotsDofs} illustrates the same story as observed in~\cite{Venas2018iao}; The increased continuity of the basis functions offered by IGA play a crucial role for improving the accuracy of the numerical solution.
To obtain higher accuracy using PML it is however important to know that the PML-thickness influences the accuracy as illustrated in~\Cref{Fig2:EnergyErrorPlotsDofs_2}. This in turn increases dofs used in order to maintain the aspect ratio of the elements. This problem is here not prevalent for the IEM~\cite{Venas2018iao}. Before pollution from the PML-thickness becomes dominant we see that the PML approximation yields solution close to the best approximation.
Finally, a comparison between two stretching functions are made in~\Cref{Fig2:EnergyErrorPlotsDofs_3}. The stretching function in~\Cref{Eq:sigmaType3} gives much better results than that of~\Cref{Eq:sigmaType2}. The PML simulation with the doubled PML thickness again follows closely the best approximation (here in the $L^2(\Gamma)$-norm).

\subsection{Cylinder}
In~\Cref{Fig:Cylinder_sweep} we motivate the implementation of the Combined Helmholtz Integral Formulation~\cite{Wu1991awr} (CHIEF) on the Collocation Conventional Boundary Integral Equation (CCBIE)~\cite{Venas2020ibe} formulation (named CCBIEC) which will be used to make a reference solution for this section.

Consider a cylinder\footnote{Note that the experimental simulation herein rotates the cylinder to be aligned with the $x$-axis. The eigenfrequencies are not altered by this transformation.} of length $L=\PI\,\si{m}$ and radius $R=\SI{1}{m}$ centered at the origin with domain
\begin{equation*}
    \Omega^-=\left\{\vec{x}\in\R\st x_1^2+x_2^2\leq R^2\quad\text{and}\quad-\frac{L}{2}\leq x_3\leq\frac{L}{2}\right\}.
\end{equation*}
\begin{figure}
	\centering
	\includegraphics[width=0.6\textwidth]{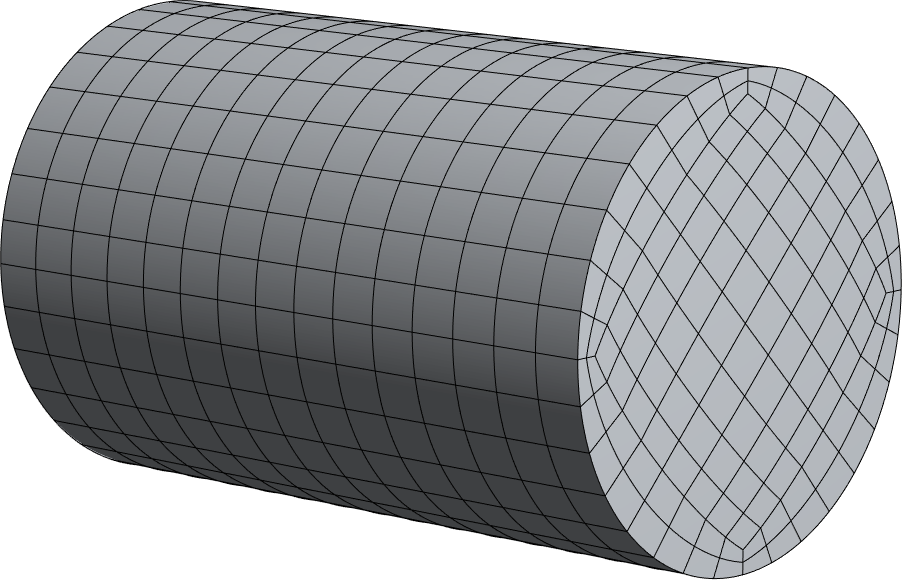}
	\caption{\textbf{Cylinder}: Boundary element method mesh ${\cal M}_{4,2,1}^{\textsc{igabem}}$ (with 736 elements and 956 degrees of freedom).}
	\label{Fig:CylinderMesh}
\end{figure}
\begin{figure}
	\centering
	\begin{subfigure}{\textwidth}
		\centering
		\includegraphics[width=0.7\textwidth]{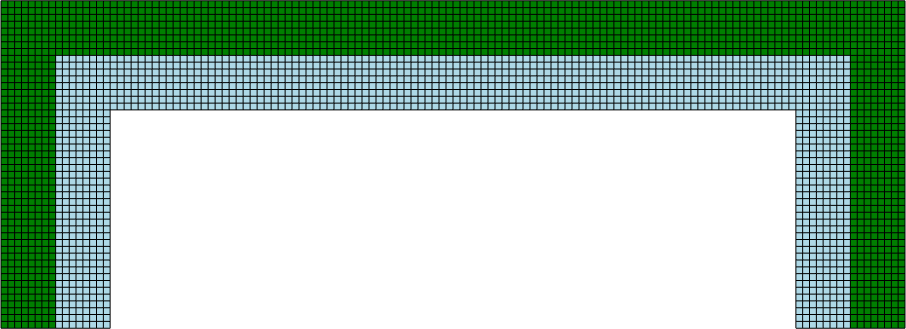}
	    \label{Fig:Cylinder_2Dmesh_pmlFill0}
		\caption{A $xz$-cross sectional of mesh ${\cal M}_{6,2,1}^{\textsc{pml}}$ (with \num{912576} elements and \num{1058544} degrees of freedom).}
    \end{subfigure}
    \par\bigskip
	\begin{subfigure}{\textwidth}
		\centering
		\includegraphics[width=0.7\textwidth]{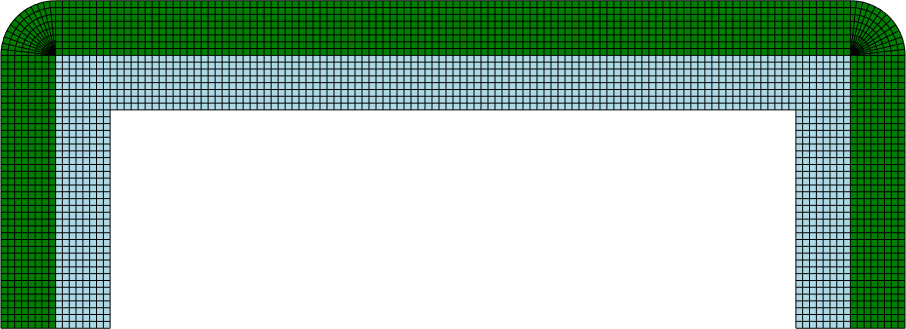}
	    \label{Fig:Cylinder_2Dmesh_pmlFill1}
		\caption{A $xz$-cross sectional of mesh ${\tilde{\cal M}}_{6,\check{p},\check{k}}^{\textsc{pml}}$ (with \num{936096} elements and \num{1082304} degrees of freedom).}
    \end{subfigure}
	\caption{\textbf{Cylinder}: Two cross sectional meshes that can be revolved around the $x$-axis (tensorially with a NURBS representation of a circle) to obtain the volumetric mesh used with the PML approach (i.e.~\Cref{Fig:CylinderMesh_PML}). The PML domain is here highlighted in green. The meshes are constructed in a way that optimizes the aspect ratio.}
	\label{Fig:Cylinder_2Dmesh}
\end{figure}
\begin{figure}
	\centering
	\includegraphics[width=0.6\textwidth]{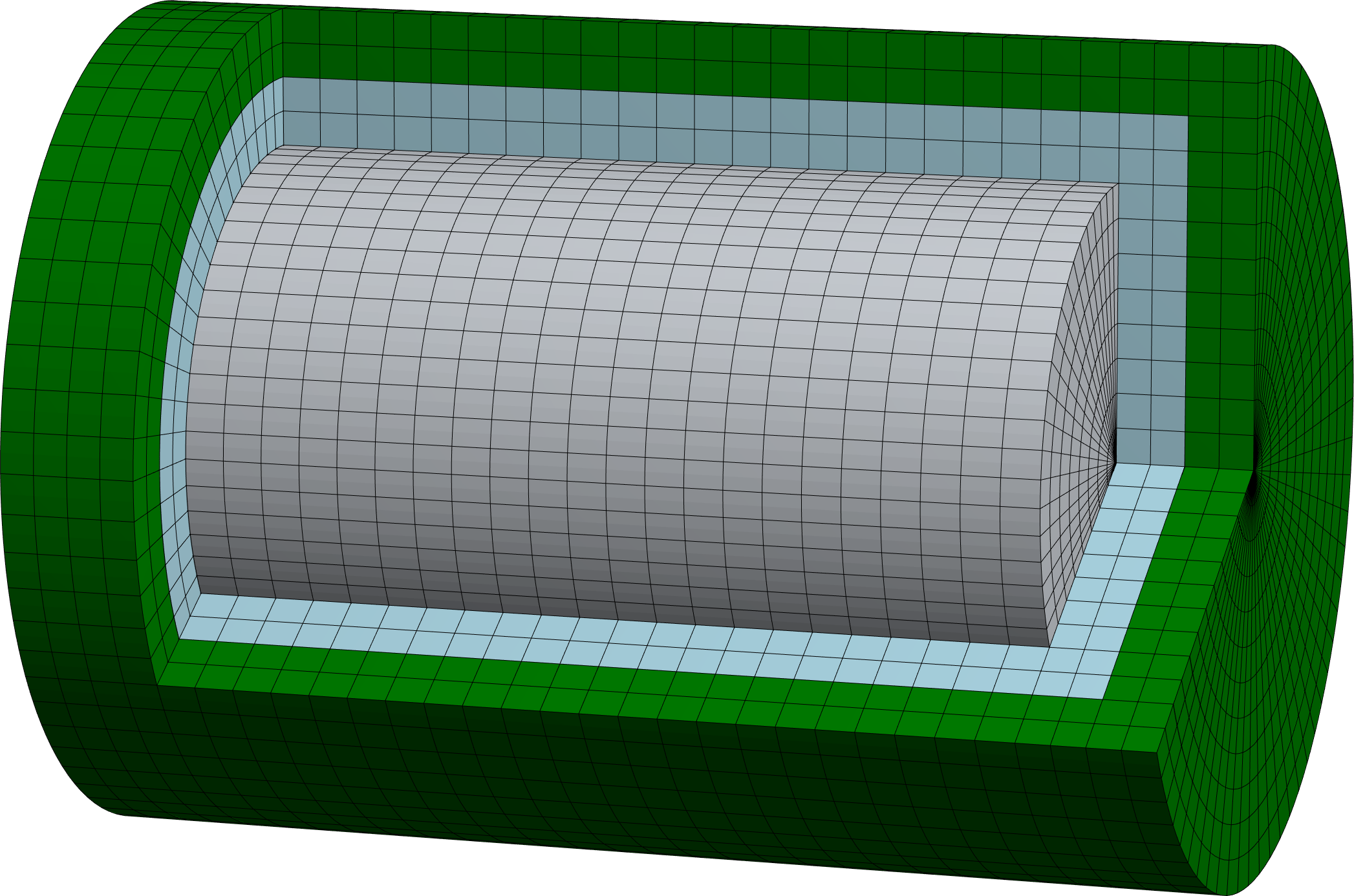}
	\caption{\textbf{Cylinder}: PML mesh ${\cal M}_{4,2,1}^{\textsc{igapml}}$ (with 12144 elements and 20712 degrees of freedom).}
	\label{Fig:CylinderMesh_PML}
\end{figure}
We want to find the eigenfunctions, $p(\vec{x})$, and eigenvalues $k^2$ to the eigenvalue problem $-\nabla^2 p = k^2 p$ with associated boundary conditions. Its interior Dirichlet problem ($p(\vec{x})=0$, $\vec{x}\in\partial\Omega^-$) has eigenfunctions\footnote{Here, $\besselJ_n$ is the cylindrical Bessel function of the first type.} (cf. \cite[p. 52]{Schenck1968iif})
\begin{equation*}
	p(\vec{x}) = \sin\frac{n_3\PI (x_3+L/2)}{L}\besselJ_n\left(\frac{x_{nm}^*}{R}r\right)\cos n\theta,\quad n_3\in\N^*,\quad \besselJ_n(x_{nm}^*)=0,\quad\vec{x}\in\Omega^-
\end{equation*}
with their corresponding wavenumbers (obtained by inserting the eigenfunctions in~\Cref{Eq:HelmholtzEqn})
\begin{equation*}
	k = \sqrt{\left(\frac{x_{nm}^*}{R}\right)^2 + \left(\frac{n_3\PI}{L}\right)^2},
\end{equation*}
where $x_{nm}^*$ excludes the trivial solutions ($x_{nm}^*\neq 0$) and the interior Neumann problem ($\partial_n p(\vec{x})=0$, $\vec{x}\in\partial\Omega^-$) has eigenfunctions
\begin{equation*}
	p(\vec{x}) = \cos\frac{n_3\PI (x_3+L/2)}{L}\besselJ_n\left(\frac{x_{nm}}{R}r\right)\cos n\theta,\quad n_3\in\N,\quad \besselJ_n'(x_{nm})=0,\quad\vec{x}\in\Omega^-
\end{equation*}
with their corresponding wavenumbers (obtained by inserting the eigenfunctions in~\Cref{Eq:HelmholtzEqn})
\begin{equation*}
	k = \sqrt{\left(\frac{x_{nm}}{R}\right)^2 + \left(\frac{n_3\PI}{L}\right)^2}.
\end{equation*}
For the exterior problem these eigenfrequencies correspond to the fictitious eigenfrequencies for the CBIE (conventional boundary integral equation) formulation and the HBIE (hyper-singular boundary integral equation) formulation, respectively. The fictitious eigenfrequencies below $kL=15$ are given by
\begin{equation}\label{Eq:CBIEfreqs}
    kL \approx 8.182137, 9.826300, 12.079081, 12.440854, 13.578794, 14.662586
\end{equation}
for the  CBIE formulation, and 
\begin{equation}
\begin{aligned}\label{Eq:HBIEfreqs}
kL &\approx 0, \PI, 5.784249, 2\PI, 6.582336, 8.540255, 3\PI, 9.595168, 10.096379, 11.058209, 11.469336, \\
&\qquad 12.037659, 12.440854, 4\PI, 13.198424, 13.449673, 13.567167, 13.578794, 13.833698, 14.617689
\end{aligned}
\end{equation}
for the HBIE formulation.

Consider the manufactured solution (cf.~\cite{Venas2020ibe})
\begin{equation}\label{Eq3:manuMFS}
	p(\vec{x}) = \sum_{n=1}^N C_n \Phi_k(\vec{x},\vec{y}_n),\quad \Phi_k(\vec{x},\vec{y}) = \frac{\euler^{\imag kR}}{4\PI R},\quad\text{where}\quad R = |\vec{x} - \vec{y}|\quad\text{and}\quad C_n = \cos(n-1).
\end{equation}
with $N=3^3=27$ source points 
\begin{equation*}
	\vec{y}_n = \frac{R}{4}[c_i, c_j, c_l],\quad n=i+3(j-1)+3^2(l-1),\quad i,j,l=1,2,3
\end{equation*}
where $c_1=-1$, $c_2=0$ and $c_3=1$. A total of $4^3=64$ uniformly spaced points around the origin on a regular cube grid of side length 1/4 is used for the interior points in the CHIEF formulation. The results are given in~\Cref{Fig:Cylinder_sweep} where we can see that the CCBIEC formulation follows the best approximation (BA) throughout the frequency sweep without any fictitious eigenfrequencies present in the CCBIE and CHBIE formulations. The CBM (collocation Burton-Miller) formulation removes the fictitious frequencies but has a reduces accuracy compared to the CCBIEC formulation.
\begin{figure}
	\centering
	\includegraphics[width=\textwidth]{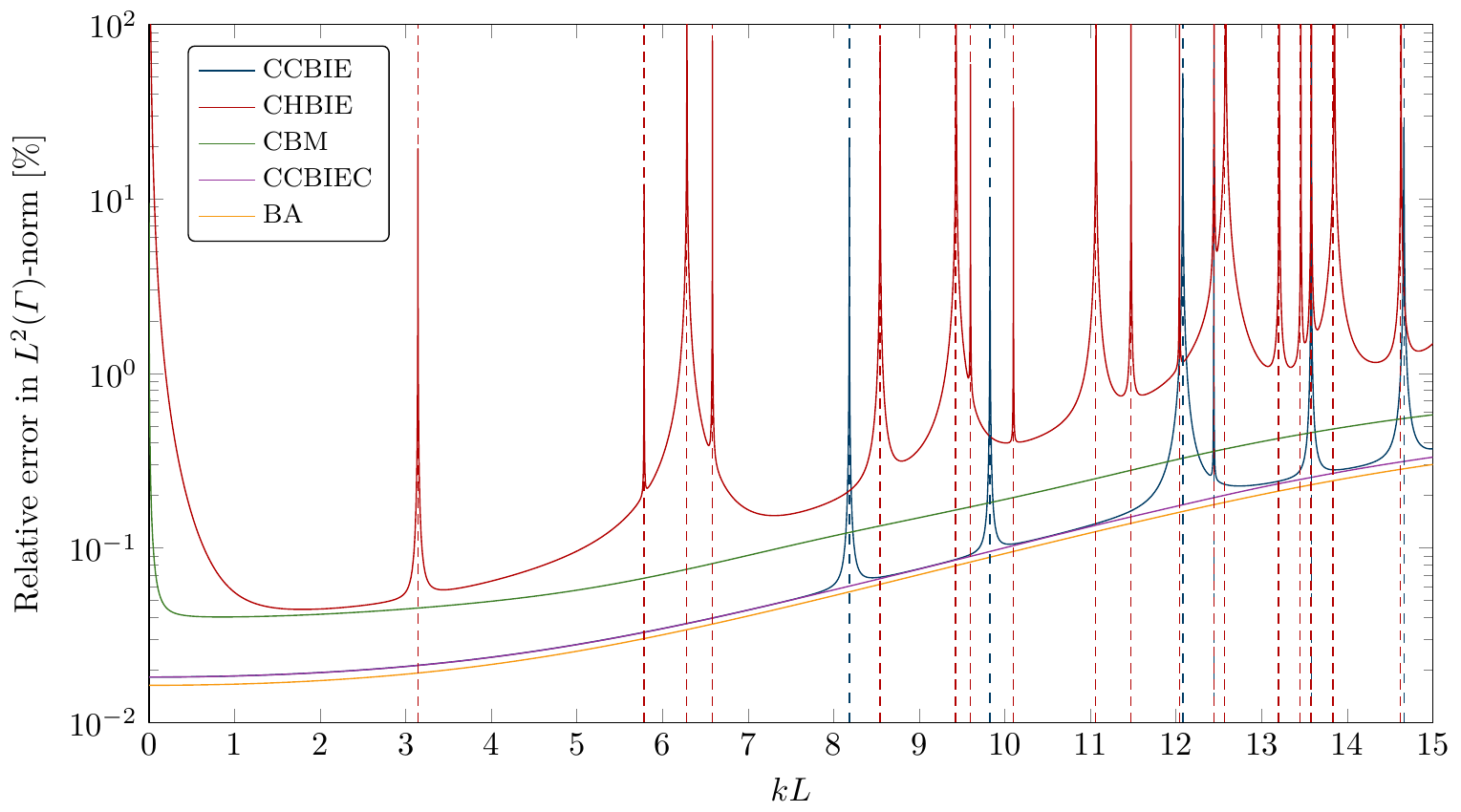}
	\caption{\textbf{Manufactured solution with a cylinder}: The plots show the instabilities around eigenfrequencies of the corresponding interior Dirichlet problem using IGABEM formulations~\cite{Venas2020ibe}. All computations are done using the parametrization in \Cref{Fig:CylinderMesh} with NURBS degree 2. The dashed lines correspond to the fictitious eigenfrequencies in~\Cref{Eq:CBIEfreqs,Eq:HBIEfreqs}.}
	\label{Fig:Cylinder_sweep}
\end{figure}
Consider now the same cylinder scattering a plane wave incident with the $x$-direction. Some PML meshes are found in~\Cref{Fig:Cylinder_2Dmesh} where two meshing strategies are outlined. Both strategies yields roughly the same accuracy. The possibility to fill corners like this will be important to create an automatic PML mesh generator for non-smooth artificial boundary as in this example. The near-field is plotted in~\Cref{Fig:Cylinder} (for $k=\SI{100}{m^{-1}}$) and the far field (for $k=\SI{50}{m^{-1}}$) in~\Cref{Fig:Cylinder_BI}. 
\begin{figure}
	\centering
	\includegraphics{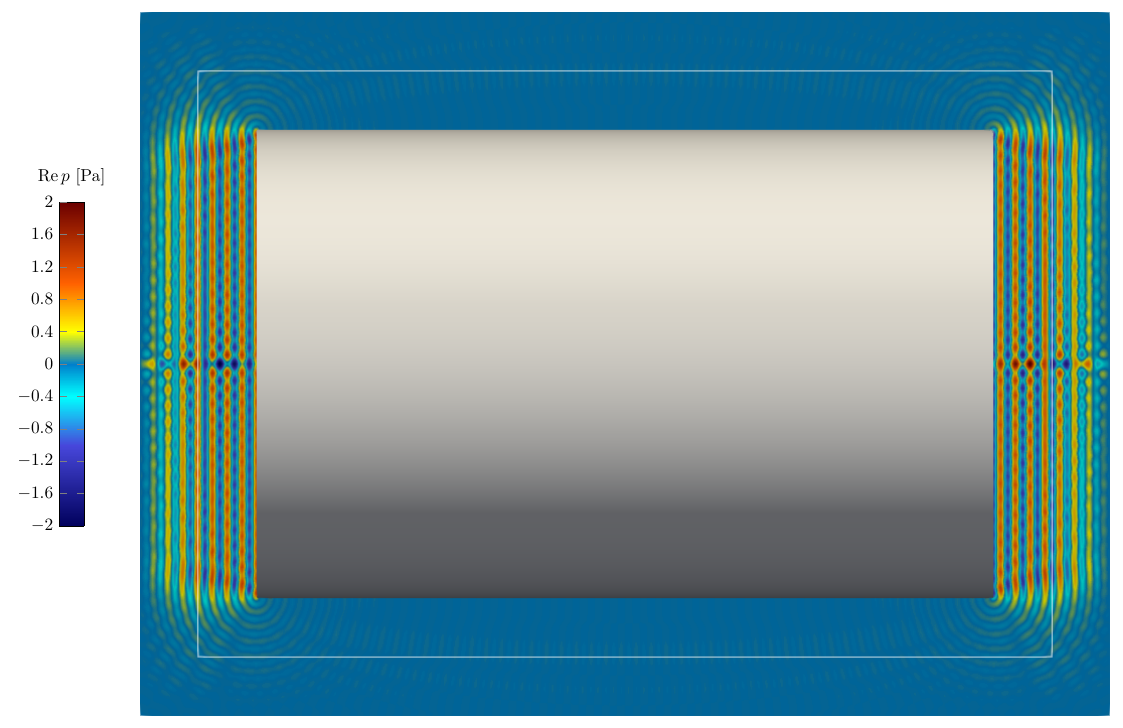}
	\caption{\textbf{Scattering from cylinder}: The figure show the scattered near field at $k=\SI{100}{m^{-1}}$ using mesh ${\cal M}_{7,2,1}^{\textsc{igapml}}$. The boundary $\Gamma_{\mathrm{a}}$ is highlighted with a white line.}
	\label{Fig:Cylinder}
\end{figure}
\begin{figure}
	\centering
	\includegraphics[width=\textwidth]{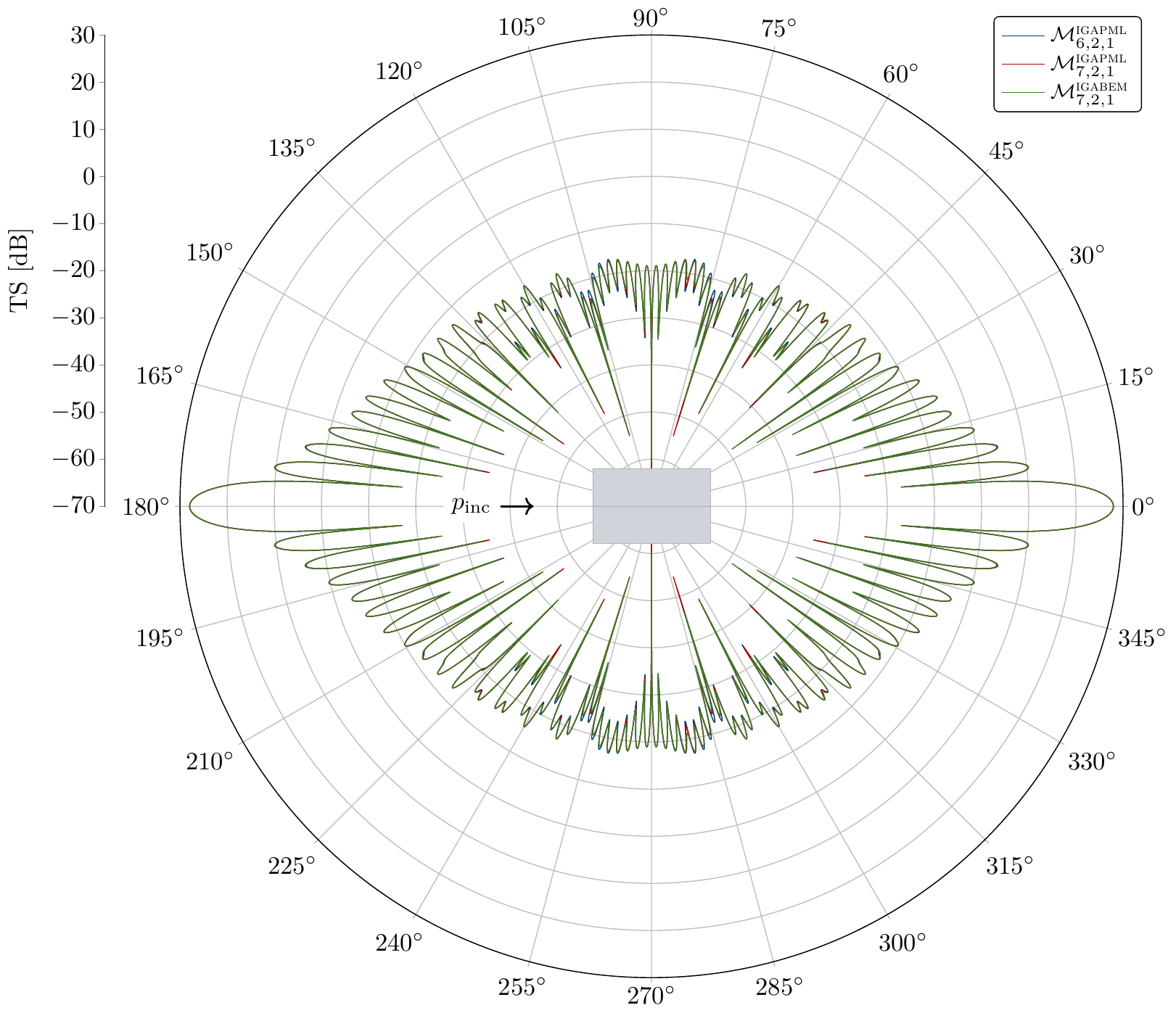}
	\caption{\textbf{Scattering from cylinder}: Far field pattern as a function of the aspect angle $\alpha$ at $k=\SI{50}{m^{-1}}$. The CCBIEC formulation was used for the BEM reference solution. The back scattering (at $\alpha=\ang{180}$) and forward scattering (at $\alpha=\ang{0}$) is present as expected.}
	\label{Fig:Cylinder_BI}
\end{figure} 

\subsection{BeTSSi Model 3}
\begin{figure}
	\centering
	\includegraphics[width=\textwidth]{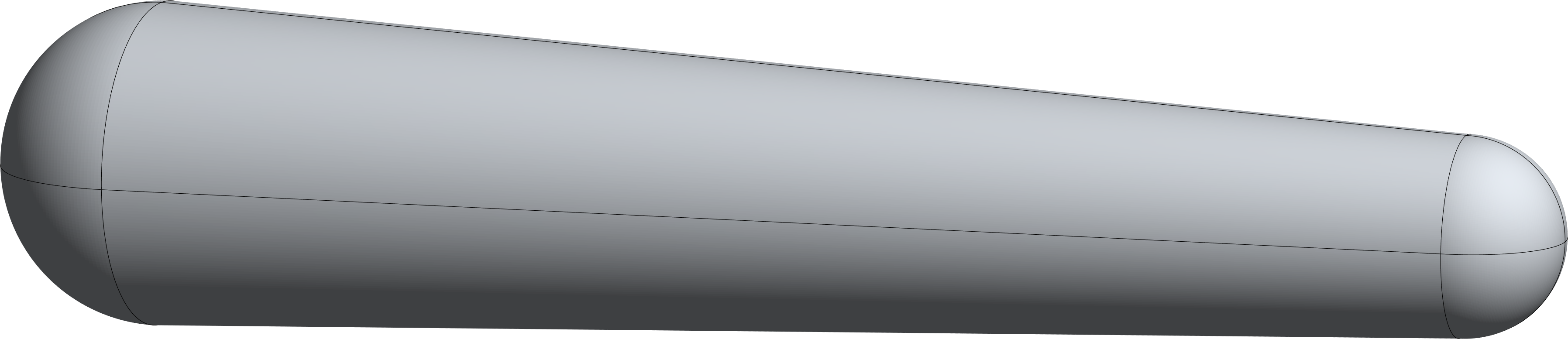}
	\caption{\textbf{The BeTSSi model 3}: This model is one of the benchmark models in the BeTSSi (Benchmark Target Strength Simulations \cite{Nolte2014bib}) community. The BeTSSi model 3 (M3) is a model given by two hemispherical end caps with radii $R_1=\SI{3}{m}$ and $R_2=\SI{5}{m}$ connected by a cone of length $L=\SI{41}{m}$. The speed of sound in the fluid is $c=\SI{1500}{ms^{-1}}$.}
	\label{Fig:M3}
\end{figure}

\begin{figure}
	\centering
	\includegraphics[width=\textwidth]{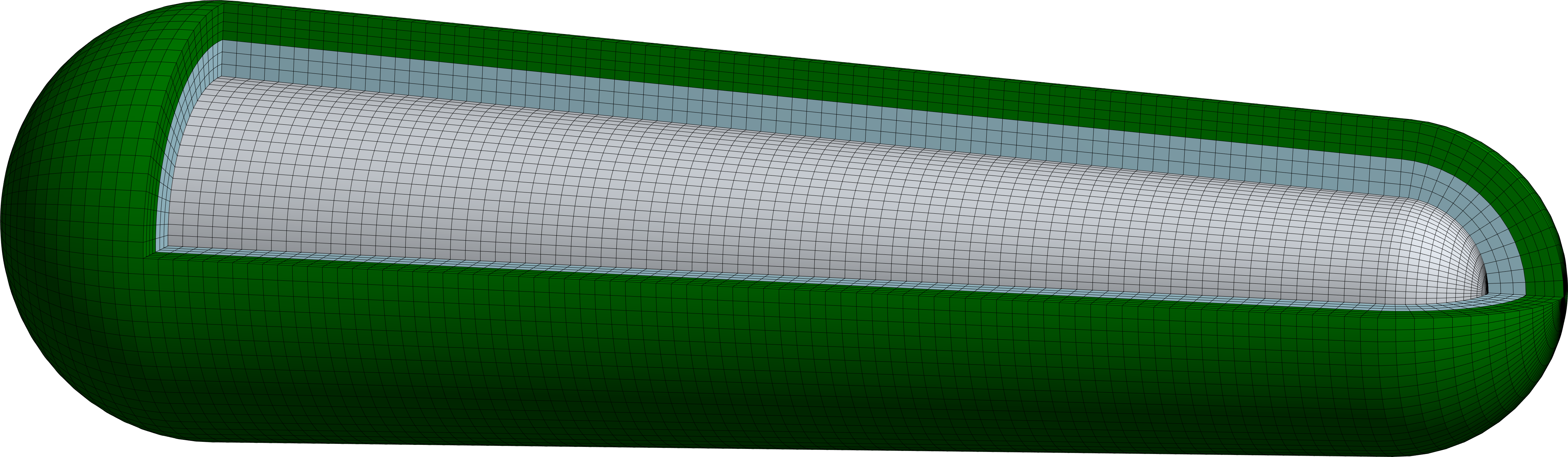}
	\caption{\textbf{The BeTSSi model 3}: Mesh ${\cal M}_{5,\check{p},\check{k}}^{\textsc{igapml}}$ illustrating the perfectly matched layer (green) domain $\Omega_{\mathrm{b}}$ around the (light blue) domain $\Omega_{\mathrm{a}}$. The distance to the PML layer is the same as the thickness; $t_{\mathrm{pml}} = 0.25 R_2$.}
	\label{Fig:M3_3Dmesh}
\end{figure}

\begin{figure}
	\centering
	\begin{subfigure}{\textwidth}
		\centering
		\includegraphics[width=\textwidth]{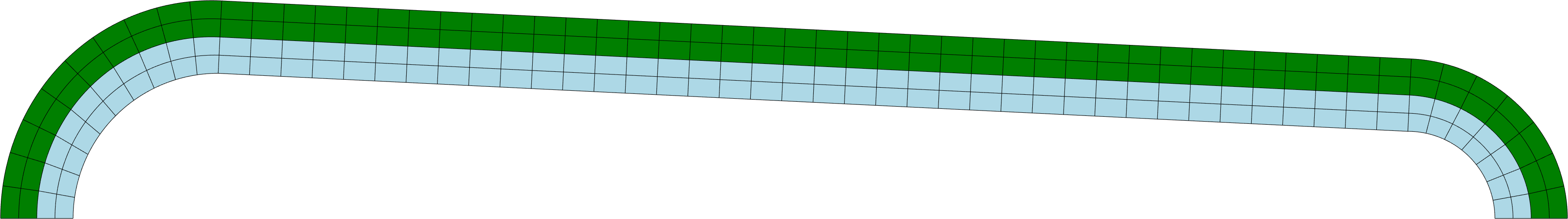}
		\caption{Mesh ${\cal M}_{4,2,1}^{\textsc{igapml}}$ with $\num{8800}$ elements and $\num{15060}$ dofs}
    \end{subfigure}
	\begin{subfigure}{\textwidth}
		\centering
		\includegraphics[width=\textwidth]{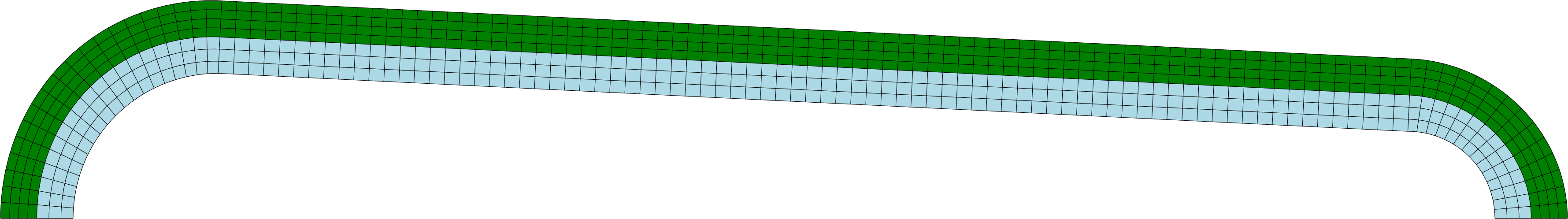}
		\caption{Mesh ${\cal M}_{5,2,1}^{\textsc{igapml}}$ with $\num{63280}$ elements and $\num{86958}$ dofs}
    \end{subfigure}
	\begin{subfigure}{\textwidth}
		\centering
		\includegraphics[width=\textwidth]{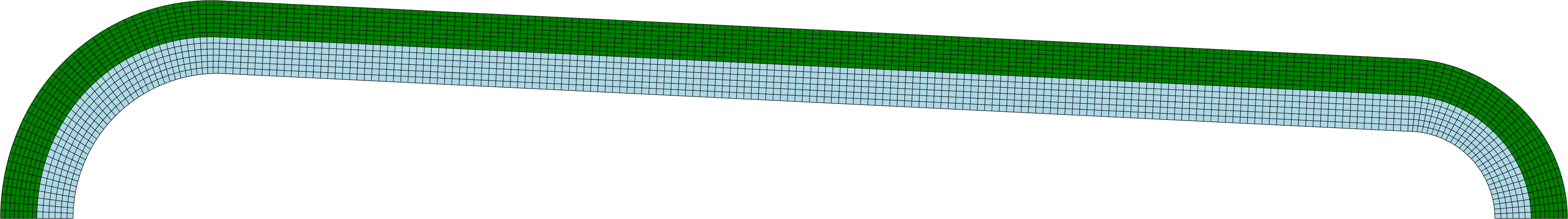}
		\caption{Mesh ${\cal M}_{6,2,1}^{\textsc{igapml}}$ with $\num{515200}$ elements and $\num{608800}$ dofs}
    \end{subfigure}
	\caption{\textbf{The BeTSSi model 3}: Illustration of the meshes where the refinement is performed reducing the aspect ratio of the coarse mesh.}
	\label{Fig:M3_2Dmesh}
\end{figure}

\begin{figure}
	\centering
	\includegraphics[width=\textwidth]{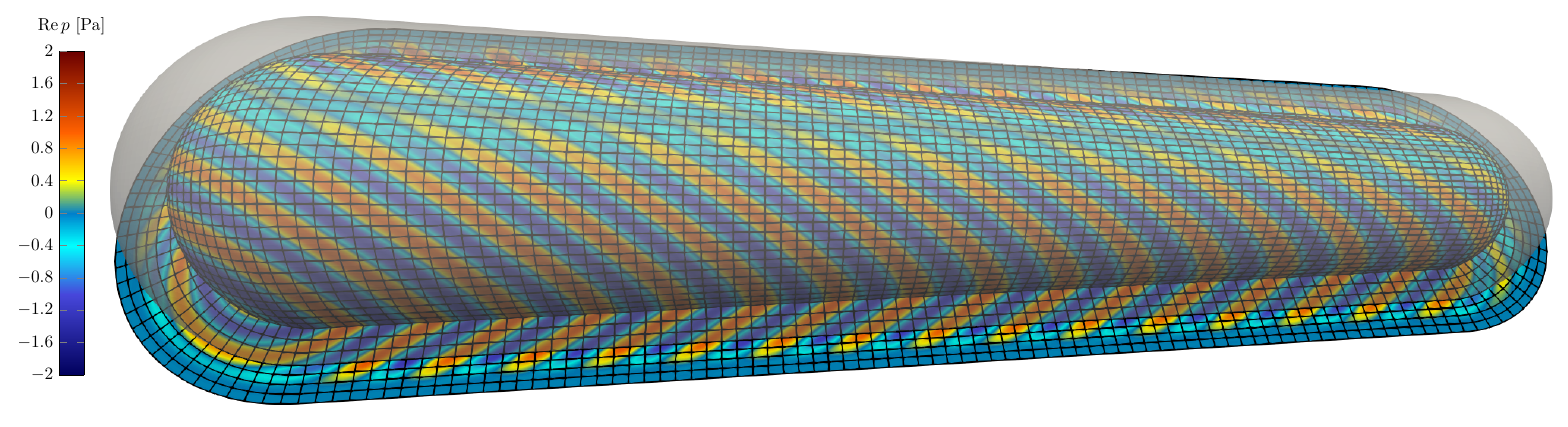}
	\caption{\textbf{The BeTSSi model 3}: Results on mesh ${\cal M}_{5,2,1}^{\textsc{igapml}}$ illustrating the scattered pressure both in the vicinity of the scatterer and in the PML. The artificial boundary $\Gamma_{\mathrm{a}}$ is added as a transparent surface.}
	\label{Fig:M3_paraviewResults}
\end{figure}

\begin{figure}
	\centering
	\includegraphics[width=\textwidth]{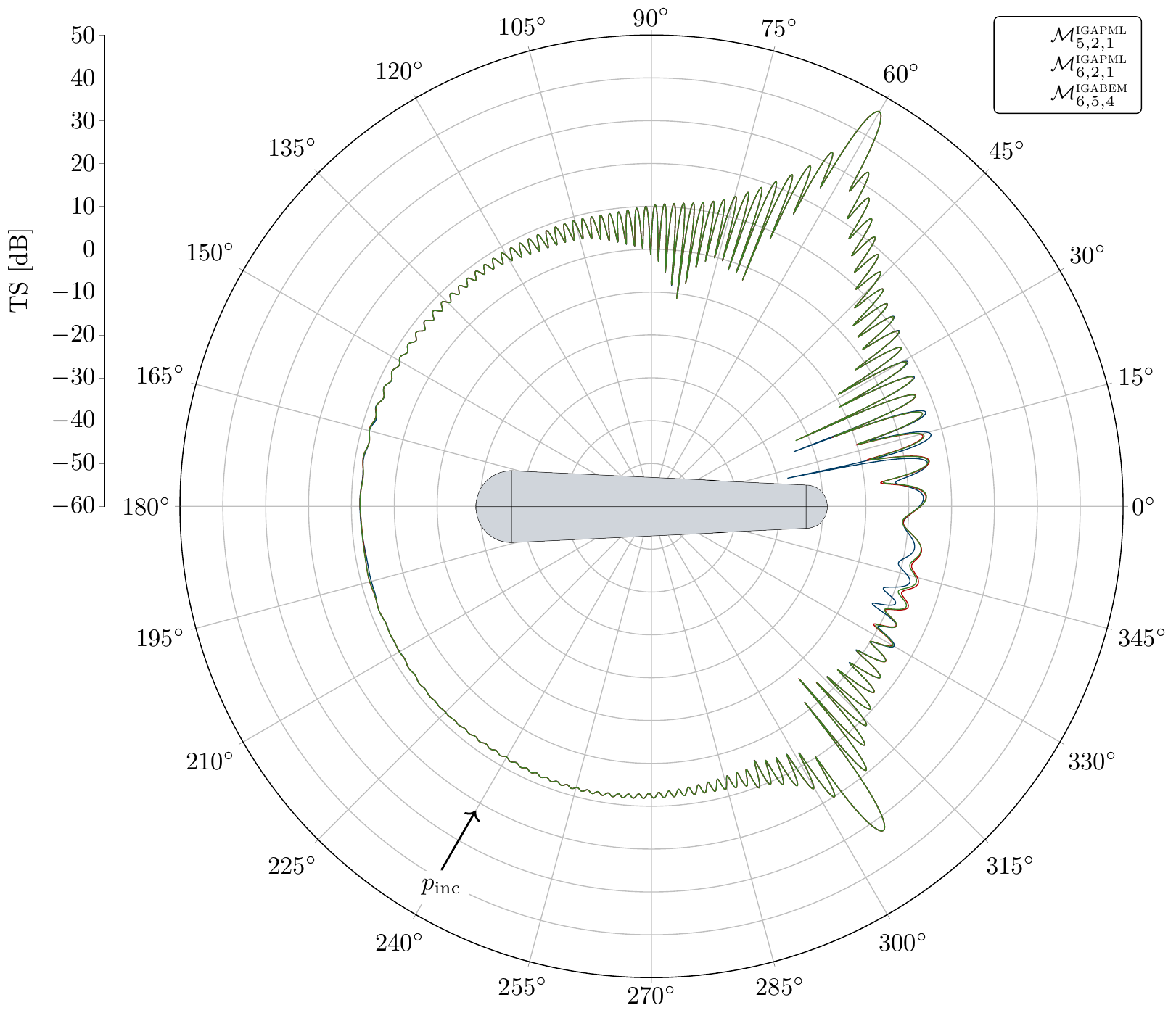}
	\caption{\textbf{The BeTSSi model 3}: Far field pattern as a function of the aspect angle $\alpha$. The CCBIE formulation was used for the BEM reference solution.}
	\label{Fig:M3_BI}
\end{figure}

\begin{figure}
	\centering
	\includegraphics[width=\textwidth]{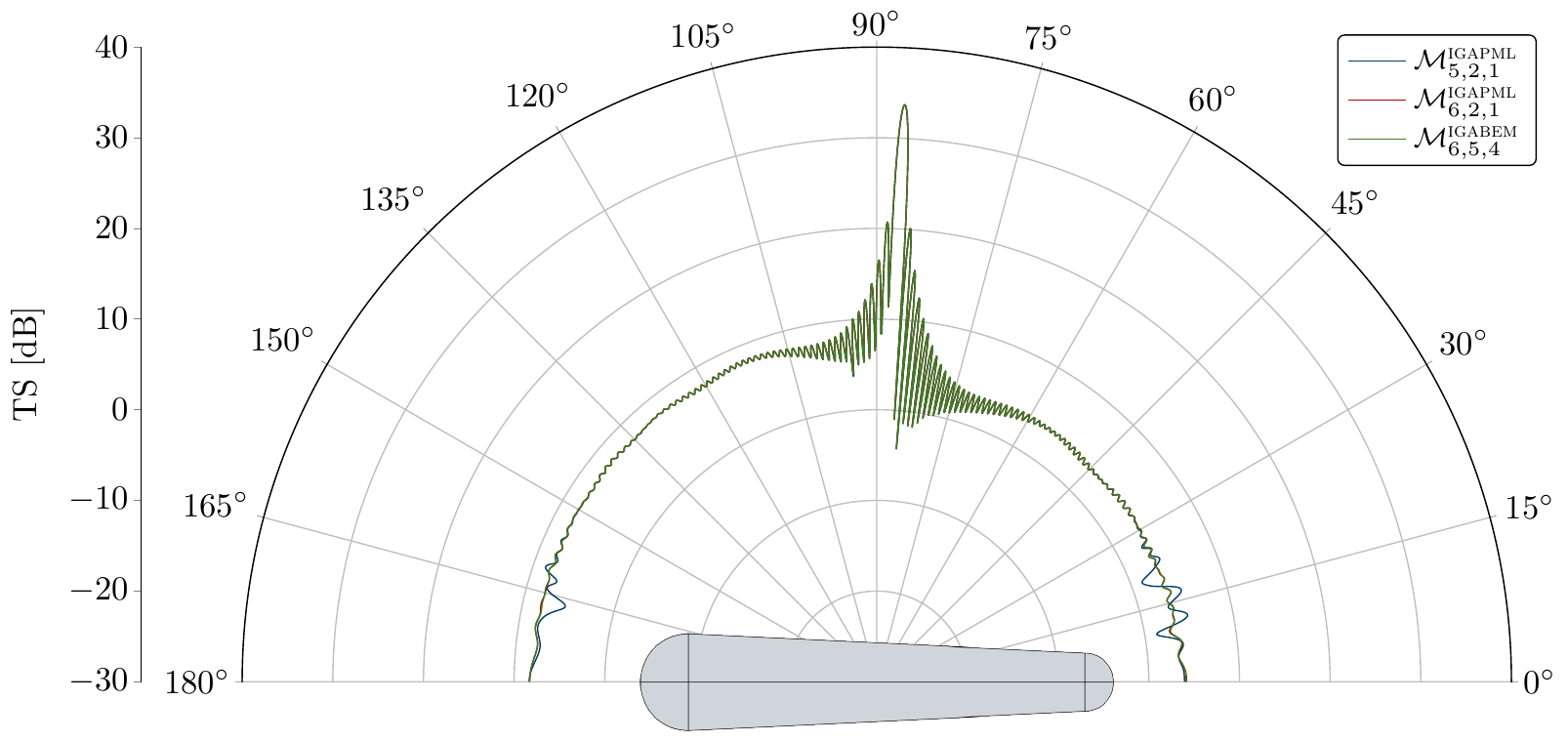}
	\caption{\textbf{The BeTSSi model 3}: Far field pattern as a function of the aspect angle $\alpha$ at $f=\SI{1}{kHz}$. }
	\label{Fig:M3_MS}
\end{figure}

The BeTSSi model 3 is illustrated in~\Cref{Fig:M3}. The model is strictly speaking neither convex nor smooth ($G^1$) due to the attachment of the smaller hemispherical cap to the intermediate cone shape, so this has been taken into account when constructing the artificial boundary.
Some meshes are visualized in~\Cref{Fig:M3_3Dmesh,Fig:M3_2Dmesh} (with polynomial order $\check{p}$ and continuity $\check{k}$), and a result for the case of an $\alpha_{\mathrm{s}}=\ang{240}$ and $\beta_{\mathrm{s}}=\ang{0}$ angle of incidence is illustrated in~\Cref{Fig:M3_paraviewResults}. 
The far field of the same simulation is plotted in~\Cref{Fig:M3_BI} as a function of the aspect angle $\alpha$. Mesh ${\cal M}_{5,2,1}^{\textsc{igapml}}$ is in good agreement with the reference solution (using BEM~\cite{Venas2020ibe}). Mesh ${\cal M}_{6,2,1}^{\textsc{igapml}}$ is visually indistinguishable with the reference solution.
Similar results are shown in a monostatic case (with $\alpha_{\mathrm{s}}=\alpha$ and $\beta_{\mathrm{s}}=\beta=0$) in~\Cref{Fig:M3_MS}.

\clearpage
\section{Conclusions}
\label{Sec:conclusions}
In this work a general spline based PML formulation has been presented to ease the construction of the PML domain. The formulation is no longer dependent on coordinate systems like the cylindrical, ellipsoidal or the Cartesian coordinate system. Instead, it is based directly on the spline space in which the numerical solution is sought. This eliminates calls to trigonometric functions and enables the PML to be truly matched to any convex and smooth artificial boundary. 

For smooth artificial boundaries this work has presented an automatic PML-layer generation approach. The outer surface of the layer located a distance $t_{\textsc{pml}}$ normally away from the inner surface of the PML-layer is approximated using least squares. These two surfaces are then linearly lofted to get the volumetric mesh. As $n_{\mathrm{dofs}}\to\infty$ this approach converges to an exact conformal PML formulation. This approach enables a simple adjustment to the standard assembly routine and requires quantities already present in standard codes without the need to evaluate NURBS-function with complex parametric arguments. The application to non-smooth artificial boundary could in principle also be implemented in an automatic fashion based on a wedge fill outlined in this work. This is suggested as future work. 

Only a modification of the standard Jacobian is needed in the bilinear form instead of several Jacobians from additional coordinate transformations. As restrictions to the computational domains are reduced, both mesh quality and computational efficiency can be improved by the present approach. Moreover, the usage of spline basis function of higher continuity improves accuracy of both the PML-layer approximation and the numerical solution itself through the IGA framework. 

The usage of unbounded absorption function conveniently reduces the number of PML parameters to tune. In the experiments we only need to tune the distance to the artificial boundary and the thickness of the PML layer. For a reasonable distance to the artificial boundary, we have used the same distance as the PML thickness yielding reasonably good results, in which case the PML thickness is the most sensitive parameter to tune for high accuracies. Until pollution from the PML-thickness occurs, the PML approximation yields solution close to the best approximation indicating that the integrations over the unbounded absorption functions are well resolved. 

\section*{Acknowledgements}
This work was supported by SINTEF Digital.

\appendix
\section{Multiple PML absorption directions}
\label{Sec:multiplePMLdirections}
For two absorption direction (say, the $\xi_2$- and the $\xi_3$-direction) we can write
\begin{equation}\label{Eq:XtwoAbs}
 	\vec{X}(\vec{\xi}) = \vec{X}(\xi_1,0,0) + \xi_2\pderiv{\vec{X}}{\xi_2}(\xi_1,0,0) + \xi_3\pderiv{\vec{X}}{\xi_3}(\xi_1,0,0) + \xi_2\xi_3\left(\pderiv{\vec{X}}{\xi_2}(\xi_1,0,1) - \pderiv{\vec{X}}{\xi_2}(\xi_1,0,0)\right).
\end{equation}
Note that we could also have used
\begin{equation}
	\pderiv{\vec{X}}{\xi_2}(\xi_1,0,1) - \pderiv{\vec{X}}{\xi_2}(\xi_1,0,0) = \pderiv{\vec{X}}{\xi_3}(\xi_1,1,0) - \pderiv{\vec{X}}{\xi_3}(\xi_1,0,0).
\end{equation}
From~\Cref{Eq:XtwoAbs} we find
\begin{align*}
	\ppderiv{\vec{X}}{\xi_2}{\xi_3}(\vec{\xi}) = \pderiv{\vec{X}}{\xi_2}(\xi_1,0,1) - \pderiv{\vec{X}}{\xi_2}(\xi_1,0,0)\\
	\pderiv{\vec{X}}{\xi_2}(\vec{\xi}) = \pderiv{\vec{X}}{\xi_2}(\xi_1,0,1) + \xi_3\ppderiv{\vec{X}}{\xi_2}{\xi_3}(\vec{\xi})\\
	\pderiv{\vec{X}}{\xi_3}(\vec{\xi}) = \pderiv{\vec{X}}{\xi_3}(\xi_1,0,1) + \xi_2\ppderiv{\vec{X}}{\xi_2}{\xi_3}(\vec{\xi})
\end{align*}
such that we can rewrite~\Cref{Eq:XtwoAbs} as
\begin{equation}
 	\vec{X}(\vec{\xi}) = \vec{X}(\xi_1,0,0) + \xi_2\pderiv{\vec{X}}{\xi_2} + \xi_3\pderiv{\vec{X}}{\xi_3} - \xi_2\xi_3\ppderiv{\vec{X}}{\xi_2}{\xi_3},
\end{equation}
which enables us to write\footnote{We take the liberty to write $I_i(\xi_i)$ as $I_i$ and $\sigma_i(\xi_i)$ as $\sigma_i$.}
\begin{align*}
 	\tilde{\vec{x}} &= \vec{X}(\xi_1,\xi_2+\imag I_2,\xi_3+\imag I_3) \\
 	&= \vec{X}(\vec{\xi}) + \imag I_2\pderiv{\vec{X}}{\xi_2} + \imag I_3\pderiv{\vec{X}}{\xi_3} - (\imag\xi_2 I_3 + \imag\xi_3 I_2 - I_2I_3)\ppderiv{\vec{X}}{\xi_2}{\xi_3}.
\end{align*}
We can then compute
\begin{align*}
 	\pderiv{\tilde{\vec{x}}}{\vec{\xi}} = \vec{J} + &\imag \left[ I_2\ppderiv{\vec{X}}{\xi_1}{\xi_2} + I_3\ppderiv{\vec{X}}{\xi_1}{\xi_3} - \left(\xi_2 I_3 + \xi_3 I_2 + \imag I_2I_3\right)\pppderiv{\vec{X}}{\xi_1}{\xi_2}{\xi_3},\right.\\
 	&~\left.\sigma_2\pderiv{\vec{X}}{\xi_2} - \sigma_2\left(\xi_3 +\imag I_3\right)\ppderiv{\vec{X}}{\xi_2}{\xi_3},\right.\\
 	&~\left.\sigma_3\pderiv{\vec{X}}{\xi_3} - \sigma_3\left(\xi_2 + \imag I_2\right)\ppderiv{\vec{X}}{\xi_2}{\xi_3}\right]
\end{align*}
where (using the same argument as before)
\begin{equation*}
	\pppderiv{\vec{X}}{\xi_1}{\xi_2}{\xi_3} = \sum_{i_1=1}^{n_1}\sum_{i_2=1}^{n_2}\sum_{i_3=1}^{n_3} \frac{1}{R_{i_1,i_2,i_3}^2(\vec{\xi})}\pderiv{R_{i_1,i_2,i_3}}{\xi_1}\pderiv{R_{i_1,i_2,i_3}}{\xi_2}\pderiv{R_{i_1,i_2,i_3}}{\xi_3}\vec{P}_{i_1,i_2,i_3}.
\end{equation*}
Similarly for linearity in all parametric direction we can write
\begin{equation}
 	\vec{X}(\vec{\xi}) = \vec{X}(\zerovec) + \xi_1\pderiv{\vec{X}}{\xi_1} + \xi_2\pderiv{\vec{X}}{\xi_2} + \xi_3\pderiv{\vec{X}}{\xi_3} - \xi_1\xi_2\ppderiv{\vec{X}}{\xi_1}{\xi_2} - \xi_1\xi_3\ppderiv{\vec{X}}{\xi_1}{\xi_3} - \xi_2\xi_3\ppderiv{\vec{X}}{\xi_2}{\xi_3} + \xi_1\xi_2\xi_3\pppderiv{\vec{X}}{\xi_1}{\xi_2}{\xi_3},
\end{equation}
which enables us to write
\begin{align*}
 	\tilde{\vec{x}} &= \vec{X}(\xi_1+\imag I_1,\xi_2+\imag I_2,\xi_3+\imag I_3) \\
 	&= \vec{X}(\vec{\xi}) + \imag I_1\pderiv{\vec{X}}{\xi_1} + \imag I_2\pderiv{\vec{X}}{\xi_2} + \imag I_3\pderiv{\vec{X}}{\xi_3} 
 	- (\imag\xi_1 I_2 + \imag\xi_2 I_1 - I_1I_2)\ppderiv{\vec{X}}{\xi_1}{\xi_2}\\
	&\quad - (\imag\xi_1 I_3 + \imag\xi_3 I_1 - I_1I_3)\ppderiv{\vec{X}}{\xi_1}{\xi_3}
	- (\imag\xi_2 I_3 + \imag\xi_3 I_2 - I_2I_3)\ppderiv{\vec{X}}{\xi_2}{\xi_3}\\
	&\quad + [(\xi_1+\imag I_1)(\xi_2+\imag I_2)(\xi_3+\imag I_3) - \xi_1\xi_2\xi_3]\pppderiv{\vec{X}}{\xi_1}{\xi_2}{\xi_3}.
\end{align*}
We can then compute
\begin{align*}
 	\pderiv{\tilde{\vec{x}}}{\vec{\xi}} = \vec{J} + &\imag\left[
 				 \pderiv{\vec{X}}{\xi_1}  
 	- (\xi_2 + \imag I_2)\ppderiv{\vec{X}}{\xi_1}{\xi_2} - (\xi_3  + \imag I_3)\ppderiv{\vec{X}}{\xi_1}{\xi_3} +  (\xi_2+\imag I_2)(\xi_3+\imag I_3)\pppderiv{\vec{X}}{\xi_1}{\xi_2}{\xi_3},\right.\\
	&\quad \left.\pderiv{\vec{X}}{\xi_2}  
 	- (\xi_1 + \imag I_1)\ppderiv{\vec{X}}{\xi_1}{\xi_2}- (\xi_3  + \imag I_3)\ppderiv{\vec{X}}{\xi_2}{\xi_3} +  (\xi_1+\imag I_1)(\xi_3+\imag I_3)\pppderiv{\vec{X}}{\xi_1}{\xi_2}{\xi_3},\right.\\
	&\quad \left.\pderiv{\vec{X}}{\xi_3}  
 	- (\xi_2 + \imag I_2)\ppderiv{\vec{X}}{\xi_2}{\xi_3}- (\xi_1 + \imag I_1)\ppderiv{\vec{X}}{\xi_1}{\xi_3} + (\xi_1+\imag I_1)(\xi_2+\imag I_2)\pppderiv{\vec{X}}{\xi_1}{\xi_2}{\xi_3}\right]\diag{\sigma_1,\sigma_2,\sigma_3}
\end{align*}

\bibliographystyle{TK_CM}
\bibliography{references}

\begin{thebibliography}{10}
\expandafter\ifx\csname url\endcsname\relax
  \def\url#1{\texttt{#1}}\fi
\expandafter\ifx\csname urlprefix\endcsname\relax\def\urlprefix{URL }\fi
\expandafter\ifx\csname href\endcsname\relax
  \def\href#1#2{#2} \def\path#1{#1}\fi

\bibitem{Sauter2011bem}
S.~A. Sauter, C.~Schwab,
  \href{http://dx.doi.org/10.1007/978-3-540-68093-2_4}{{\em Boundary Element
  Methods}}, Springer Berlin Heidelberg, Berlin, Heidelberg, 2011, pp.
  183--287.
\newblock

\bibitem{Schanz2007bea}
M.~Schanz, O.~Steinbach,
  \href{https://books.google.no/books?id=e1J9bdgAX94C}{{\em Boundary Element
  Analysis: Mathematical Aspects and Applications}}, Lecture Notes in Applied
  and Computational Mechanics, Springer Berlin Heidelberg, 2007.

\bibitem{Marburg2008cao}
S.~Marburg, B.~Nolte, \href{https://www.doi.org/10.1007/978-3-540-77448-8}{{\em
  Computational Acoustics of Noise Propagation in Fluids-Finite and Boundary
  Element Methods}}, vol. 578, Springer, 2008.
\newblock

\bibitem{Chandler_Wilde2012nab}
S.~N. Chandler{-}Wilde, I.~G. Graham, S.~Langdon, E.~A. Spence,
  \href{https://doi.org/10.1017/s0962492912000037}{Numerical-asymptotic
  boundary integral methods in high-frequency acoustic scattering}, {\em Acta
  Numerica}, 21:89--305 (2012).
\newblock

\bibitem{Berenger1994apm}
J.-P. Berenger, \href{https://doi.org/10.1006/jcph.1994.1159}{A perfectly
  matched layer for the absorption of electromagnetic waves}, {\em Journal of
  Computational Physics}, 114:185--200 (1994).
\newblock

\bibitem{Berenger1996pml}
J.-P. Berenger, \href{https://doi.org/10.1109/8.477535}{Perfectly matched layer
  for the {FDTD} solution of wave-structure interaction problems}, {\em {IEEE}
  Transactions on Antennas and Propagation}, 44:110--117 (1996).
\newblock

\bibitem{Matuszyk2012pfe}
P.~J. Matuszyk, L.~F. Demkowicz,
  \href{https://dx.doi.org/10.1007/s00466-012-0702-1}{Parametric finite
  elements, exact sequences and perfectly matched layers}, {\em Computational
  Mechanics}, 51:35--45 (2012).
\newblock

\bibitem{Givoli2013nmf}
D.~Givoli,
  \href{https://books.google.no/books?hl=en&lr=&id=e0svBQAAQBAJ&oi=fnd&pg=PP1&dq=+D.+Givoli,+Numerical+Methods+for+Problems+in+Infinite+Domains,+El-+sevier,+Amsterdam,+1992.&ots=joKWS0lPvD&sig=ppGk7zCR_o36wfHvgQIbvj2YTfs&redir_esc=y#v=onepage&q=D.%20Givoli%2C%20Numerical%20Methods%20for%20Problems%20in%20Infinite%20Domains%2C%20El-%20sevier%2C%20Amsterdam%2C%201992.&f=false}{{\em
  Numerical methods for problems in infinite domains}}, vol.~33, Elsevier,
  2013.

\bibitem{Shirron1995soe}
J.~J. Shirron,
  \href{http://adsabs.harvard.edu/abs/1995PhDT.......276S}{Solution of exterior
  {H}elmholtz problems using finite and infinite elements}, Ph.D. thesis,
  University of Maryland College Park (1995).

\bibitem{Bayliss1982bcf}
A.~Bayliss, M.~Gunzburger, E.~Turkel,
  \href{https://doi.org/10.1137/0142032}{Boundary conditions for the numerical
  solution of elliptic equations in exterior regions}, {\em {SIAM} Journal on
  Applied Mathematics}, 42:430--451 (1982).
\newblock

\bibitem{Hagstrom1998afo}
T.~Hagstrom, S.~Hariharan,
  \href{https://doi.org/10.1016/S0168-9274(98)00022-1}{A formulation of
  asymptotic and exact boundary conditions using local operators}, {\em Applied
  Numerical Mathematics}, 27:403--416 (1998), Special Issue on Absorbing
  Boundary Conditions.
\newblock

\bibitem{Tezaur2001tdf}
R.~Tezaur, A.~Macedo, C.~Farhat, R.~Djellouli,
  \href{http://dx.doi.org/10.1002/nme.346}{Three-dimensional finite element
  calculations in acoustic scattering using arbitrarily shaped convex
  artificial boundaries}, {\em International Journal for Numerical Methods in
  Engineering}, 53:1461--1476 (2001).
\newblock

\bibitem{Bettess1977ie}
P.~Bettess, \href{http://dx.doi.org/10.1002/nme.1620110107}{Infinite elements},
  {\em International Journal for Numerical Methods in Engineering}, 11:53--64
  (1977).
\newblock

\bibitem{Bettess1977dar}
P.~Bettess, O.~C. Zienkiewicz,
  \href{https://doi.org/10.1002/nme.1620110808}{Diffraction and refraction of
  surface waves using finite and infinite elements}, {\em International Journal
  for Numerical Methods in Engineering}, 11:1271--1290 (1977).
\newblock

\bibitem{Demkowicz2001aoa}
L.~Demkowicz, F.~Ihlenburg,
  \href{http://dx.doi.org/10.1007/PL00005440}{Analysis of a coupled
  finite-infinite element method for exterior {H}elmholtz problems}, {\em
  Numerische Mathematik}, 88:43--73 (2001).
\newblock

\bibitem{Hughes2005iac}
T.~Hughes, J.~Cottrell, Y.~Bazilevs,
  \href{https://doi.org/10.1016/j.cma.2004.10.008}{Isogeometric analysis:
  {CAD}, finite elements, {NURBS}, exact geometry and mesh refinement}, {\em
  Computer Methods in Applied Mechanics and Engineering}, 194:4135--4195
  (2005).
\newblock

\bibitem{Venas2018iao}
J.~V. Ven{\aa}s, T.~Kvamsdal, T.~Jenserud,
  \href{http://www.sciencedirect.com/science/article/pii/S004578251830094X}{Isogeometric
  analysis of acoustic scattering using infinite elements}, {\em Computer
  Methods in Applied Mechanics and Engineering}, 335:152--193 (2018).
\newblock

\bibitem{Venas2020ibe}
J.~V. Ven{\aa}s, T.~Kvamsdal,
  \href{https://dx.doi.org/10.1016/j.cma.2019.112670}{Isogeometric boundary
  element method for acoustic scattering by a submarine}, {\em Computer Methods
  in Applied Mechanics and Engineering}, 359:112670 (2020).
\newblock

\bibitem{Simpson2012atd}
R.~Simpson, S.~Bordas, J.~Trevelyan, T.~Rabczuk,
  \href{https://dx.doi.org/10.1016/j.cma.2011.08.008}{A two-dimensional
  isogeometric boundary element method for elastostatic analysis}, {\em
  Computer Methods in Applied Mechanics and Engineering}, 209-212:87--100
  (2012).
\newblock

\bibitem{Simpson2014aib}
R.~N. Simpson, M.~A. Scott, M.~Taus, D.~C. Thomas, H.~Lian,
  \href{https://doi.org/10.1016/j.cma.2013.10.026}{Acoustic isogeometric
  boundary element analysis}, {\em Computer Methods in Applied Mechanics and
  Engineering}, 269:265--290 (2014).
\newblock

\bibitem{Peake2013eib}
M.~Peake, J.~Trevelyan, G.~Coates,
  \href{https://doi.org/10.1016/j.cma.2013.03.016}{Extended isogeometric
  boundary element method ({XIBEM}) for two-dimensional {H}elmholtz problems},
  {\em Computer Methods in Applied Mechanics and Engineering}, 259:93--102
  (2013).
\newblock

\bibitem{Peake2014eai}
M.~Peake, \href{http://etheses.dur.ac.uk/10655/}{Enriched and isogeometric
  boundary element methods for acoustic wave scattering}, Ph.D. thesis, Durham
  University (2014).

\bibitem{Peake2015eib}
M.~Peake, J.~Trevelyan, G.~Coates,
  \href{https://doi.org/10.1016/j.cma.2014.10.039}{Extended isogeometric
  boundary element method ({XIBEM}) for three-dimensional medium-wave acoustic
  scattering problems}, {\em Computer Methods in Applied Mechanics and
  Engineering}, 284:762--780 (2015), {Isogeometric Analysis Special Issue}.
\newblock

\bibitem{Dolz2016aib}
J.~D\"{o}lz, H.~Harbrecht, M.~Peters,
  \href{https://onlinelibrary.wiley.com/doi/abs/10.1002/nme.5274}{An
  interpolation-based fast multipole method for higher-order boundary elements
  on parametric surfaces}, {\em International Journal for Numerical Methods in
  Engineering}, 108:1705--1728 (2016).
\newblock

\bibitem{Coox2017aii}
L.~Coox, O.~Atak, D.~Vandepitte, W.~Desmet,
  \href{https://doi.org/10.1016/j.cma.2016.05.039}{An isogeometric indirect
  boundary element method for solving acoustic problems in open-boundary
  domains}, {\em Computer Methods in Applied Mechanics and Engineering},
  316:186--208 (2017).
\newblock

\bibitem{Keuchel2017eoh}
S.~Keuchel, N.~C. Hagelstein, O.~Zaleski, O.~von Estorff,
  \href{http://www.sciencedirect.com/science/article/pii/S0045782517305704}{Evaluation
  of hypersingular and nearly singular integrals in the isogeometric boundary
  element method for acoustics}, {\em Computer Methods in Applied Mechanics and
  Engineering}, 325:488--504 (2017).
\newblock

\bibitem{Dolz2018afi}
J.~D\"{o}lz, H.~Harbrecht, S.~Kurz, S.~Sch\"{o}ps, F.~Wolf,
  \href{http://www.sciencedirect.com/science/article/pii/S0045782517306916}{A
  fast isogeometric {BEM} for the three dimensional {L}aplace- and {H}elmholtz
  problems}, {\em Computer Methods in Applied Mechanics and Engineering},
  330:83--101 (2018).
\newblock

\bibitem{Sun2019dib}
Y.~Sun, J.~Trevelyan, G.~Hattori, C.~Lu,
  \href{http://www.sciencedirect.com/science/article/pii/S0955799718304089}{Discontinuous
  isogeometric boundary element ({IGABEM}) formulations in 3{D} automotive
  acoustics}, {\em Engineering Analysis with Boundary Elements}, 105:303--311
  (2019).
\newblock

\bibitem{Wu2020iib}
Y.~Wu, C.~Dong, H.~Yang,
  \href{http://www.sciencedirect.com/science/article/pii/S0377042719302869}{Isogeometric
  indirect boundary element method for solving the 3{D} acoustic problems},
  {\em Journal of Computational and Applied Mathematics}, 363:273--299 (2020).
\newblock

\bibitem{Liu2017soo}
C.~Liu, L.~Chen, W.~Zhao, H.~Chen,
  \href{https://dx.doi.org/10.1016/j.enganabound.2017.09.009}{Shape
  optimization of sound barrier using an isogeometric fast multipole boundary
  element method in two dimensions}, {\em Engineering Analysis with Boundary
  Elements}, 85:142--157 (2017).
\newblock

\bibitem{Chen2018aia}
L.~Chen, C.~Liu, W.~Zhao, L.~Liu,
  \href{https://dx.doi.org/10.1016/j.cma.2018.03.025}{An isogeometric approach
  of two dimensional acoustic design sensitivity analysis and topology
  optimization analysis for absorbing material distribution}, {\em Computer
  Methods in Applied Mechanics and Engineering}, 336:507--532 (2018).
\newblock

\bibitem{Chen2019ifm}
L.~Chen, W.~Zhao, C.~Liu, H.~Chen, S.~Marburg,
  \href{http://journals.pan.pl/Content/112810/PDF/aoa.2019.129263.pdf}{Isogeometric
  fast multipole boundary element method based on burton-miller formulation for
  3d acoustic problems}, {\em Archives of Acoustics}, vol. 44:475--492 (2019).
\newblock

\bibitem{Chen2019sso}
L.~Chen, H.~Lian, Z.~Liu, H.~Chen, E.~Atroshchenko, S.~Bordas,
  \href{https://dx.doi.org/10.1016/j.cma.2019.06.012}{Structural shape
  optimization of three dimensional acoustic problems with isogeometric
  boundary element methods}, {\em Computer Methods in Applied Mechanics and
  Engineering}, 355:926--951 (2019).
\newblock

\bibitem{Chen2020ato}
L.~Chen, C.~Lu, H.~Lian, Z.~Liu, W.~Zhao, S.~Li, H.~Chen, S.~P. Bordas,
  \href{https://dx.doi.org/10.1016/j.cma.2019.112806}{Acoustic topology
  optimization of sound absorbing materials directly from subdivision surfaces
  with isogeometric boundary element methods}, {\em Computer Methods in Applied
  Mechanics and Engineering}, 362:112806 (2020).
\newblock

\bibitem{Shaaban2020sob}
A.~M. Shaaban, C.~Anitescu, E.~Atroshchenko, T.~Rabczuk,
  \href{https://dx.doi.org/10.1016/j.enganabound.2019.12.012}{Shape
  optimization by conventional and extended isogeometric boundary element
  method with {PSO} for two-dimensional helmholtz acoustic problems}, {\em
  Engineering Analysis with Boundary Elements}, 113:156--169 (2020).
\newblock

\bibitem{Shaaban2020ibe}
A.~M. Shaaban, C.~Anitescu, E.~Atroshchenko, T.~Rabczuk,
  \href{https://dx.doi.org/10.1016/j.jsv.2020.115598}{Isogeometric boundary
  element analysis and shape optimization by {PSO} for 3d axi-symmetric high
  frequency helmholtz acoustic problems}, {\em Journal of Sound and Vibration},
  486:115598 (2020).
\newblock

\bibitem{Shaaban20213ib}
A.~M. Shaaban, C.~Anitescu, E.~Atroshchenko, T.~Rabczuk,
  \href{https://dx.doi.org/10.1016/j.cma.2021.113950}{3d isogeometric boundary
  element analysis and structural shape optimization for helmholtz acoustic
  scattering problems}, {\em Computer Methods in Applied Mechanics and
  Engineering}, 384:113950 (2021).
\newblock

\bibitem{Shaaban2022aib}
A.~M. Shaaban, C.~Anitescu, E.~Atroshchenko, T.~Rabczuk,
  \href{https://dx.doi.org/10.1016/j.apacoust.2021.108410}{An isogeometric
  burton-miller method for the transmission loss optimization with application
  to mufflers with internal extended tubes}, {\em Applied Acoustics},
  185:108410 (2022).
\newblock

\bibitem{Xie2021aam}
X.~Xie, Y.~Liu, \href{https://dx.doi.org/10.1016/j.cma.2020.113532}{An adaptive
  model order reduction method for boundary element-based multi-frequency
  acoustic wave problems}, {\em Computer Methods in Applied Mechanics and
  Engineering}, 373:113532 (2021).
\newblock

\bibitem{Hetmaniuk2012raa}
U.~Hetmaniuk, R.~Tezaur, C.~Farhat,
  \href{https://onlinelibrary.wiley.com/doi/abs/10.1002/nme.4271}{Review and
  assessment of interpolatory model order reduction methods for frequency
  response structural dynamics and acoustics problems}, {\em International
  Journal for Numerical Methods in Engineering}, 90:1636--1662 (2012).
\newblock

\bibitem{Hetmaniuk2013aas}
U.~Hetmaniuk, R.~Tezaur, C.~Farhat,
  \href{https://onlinelibrary.wiley.com/doi/abs/10.1002/nme.4436}{An adaptive
  scheme for a class of interpolatory model reduction methods for frequency
  response problems}, {\em International Journal for Numerical Methods in
  Engineering}, 93:1109--1124 (2012).
\newblock

\bibitem{Quarteroni2015rbm}
A.~Quarteroni, A.~Manzoni, F.~Negri,
  \href{https://books.google.no/books?id=e6FnCgAAQBAJ}{{\em Reduced Basis
  Methods for Partial Differential Equations: An Introduction}}, UNITEXT,
  Springer International Publishing, 2015.

\bibitem{Safjan2002tdi}
A.~Safjan, M.~Newman,
  \href{http://www.sciencedirect.com/science/article/pii/S0898122102800086}{Three-dimensional
  infinite elements utilizing basis functions with compact support}, {\em
  Computers {\&} Mathematics with Applications}, 43:981--1002 (2002).
\newblock

\bibitem{Beriot2020aap}
H.~B{\'{e}}riot, A.~Modave, \href{https://dx.doi.org/10.1002/nme.6560}{An
  automatic perfectly matched layer for acoustic finite element simulations in
  convex domains of general shape}, {\em International Journal for Numerical
  Methods in Engineering} (2020).
\newblock

\bibitem{Ozgun2007nml}
O.~Ozgun, M.~Kuzuoglu,
  \href{https://dx.doi.org/10.1109/tap.2007.891865}{Non-maxwellian
  locally-conformal {PML} absorbers for finite element mesh truncation}, {\em
  {IEEE} Transactions on Antennas and Propagation}, 55:931--937 (2007).
\newblock

\bibitem{Ozgun2007nfp}
O.~Ozgun, M.~Kuzuoglu,
  \href{https://www.sciencedirect.com/science/article/pii/S0021999107003919}{Near-field
  performance analysis of locally-conformal perfectly matched absorbers via
  monte carlo simulations}, {\em Journal of Computational Physics},
  227:1225--1245 (2007).
\newblock

\bibitem{Mi2021ilc}
Y.~Mi, X.~Yu, \href{https://dx.doi.org/10.1016/j.cma.2021.113925}{Isogeometric
  locally-conformal perfectly matched layer for time-harmonic acoustics}, {\em
  Computer Methods in Applied Mechanics and Engineering}, 384:113925 (2021).
\newblock

\bibitem{Drzisga2020tsm}
D.~Drzisga, B.~Keith, B.~Wohlmuth,
  \href{https://dx.doi.org/10.1016/j.cma.2020.113322}{The surrogate matrix
  methodology: Accelerating isogeometric analysis of waves}, {\em Computer
  Methods in Applied Mechanics and Engineering}, 372:113322 (2020).
\newblock

\bibitem{Marussig2017aro}
B.~Marussig, T.~J.~R. Hughes,
  \href{https://dx.doi.org/10.1007/s11831-017-9220-9}{A review of trimming in
  isogeometric analysis: Challenges, data exchange and simulation aspects},
  {\em Archives of Computational Methods in Engineering}, 25:1059--1127 (2017).
\newblock

\bibitem{Urick2019wbo}
B.~Urick, B.~Marussig, E.~Cohen, R.~H. Crawford, T.~J. Hughes, R.~F.
  Riesenfeld, \href{https://dx.doi.org/10.1016/j.cad.2019.05.034}{Watertight
  boolean operations: A framework for creating {CAD}-compatible gap-free
  editable solid models}, {\em Computer-Aided Design}, 115:147--160 (2019).
\newblock

\bibitem{Hiemstra2020tun}
R.~R. Hiemstra, K.~M. Shepherd, M.~J. Johnson, L.~Quan, T.~J. Hughes,
  \href{https://dx.doi.org/10.1016/j.cma.2020.113227}{Towards untrimmed
  {NURBS}: {CAD} embedded reparameterization of trimmed b-rep geometry using
  frame-field~guided global parameterization}, {\em Computer Methods in Applied
  Mechanics and Engineering}, 369:113227 (2020).
\newblock

\bibitem{Sommerfeld1949pde}
A.~Sommerfeld,
  \href{http://www.sciencedirect.com/science/bookseries/00798169/1}{{\em
  Partial differential equations in physics}}, vol.~1, Academic press, 1949.

\bibitem{Ihlenburg1998fea}
F.~Ihlenburg, \href{https://doi.org/10.1007/b98828}{{\em Finite Element
  Analysis of Acoustic Scattering}}, vol. 132 of {\em Applied Mathematical
  Sciences}, Springer, New York, USA, 1998.
\newblock

\bibitem{Cottrell2006iao}
J.~Cottrell, A.~Reali, Y.~Bazilevs, T.~Hughes,
  \href{https://doi.org/10.1016/j.cma.2005.09.027}{Isogeometric analysis of
  structural vibrations}, {\em Computer Methods in Applied Mechanics and
  Engineering}, 195:5257--5296 (2006).
\newblock

\bibitem{Venas2015iao}
J.~V. Ven{\aa}s, \href{http://hdl.handle.net/11250/2352617}{Isogeometric
  analysis of acoustic scattering}, Master's thesis, Norwegian University of
  Science and Technology, {Trondheim, Norway} (2015).

\bibitem{Shirron2006afe}
J.~J. Shirron, T.~E. Giddings,
  \href{https://dx.doi.org/10.1016/j.cma.2006.07.009}{A finite element model
  for acoustic scattering from objects near a fluid{\textendash}fluid
  interface}, {\em Computer Methods in Applied Mechanics and Engineering},
  196:279--288 (2006).
\newblock

\bibitem{Michler2007itp}
C.~Michler, L.~Demkowicz, J.~Kurtz, D.~Pardo,
  \href{https://doi.org/10.1002/num.20252}{Improving the performance of
  perfectly matched layers by means of $hp$-adaptivity}, {\em Numerical Methods
  for Partial Differential Equations}, 23:832--858 (2007).
\newblock

\bibitem{Astaneh2018opm}
A.~V. Astaneh, B.~Keith, L.~Demkowicz,
  \href{https://dx.doi.org/10.1007/s00466-018-1640-3}{On perfectly matched
  layers for discontinuous {P}etrov{\textendash}{G}alerkin methods}, {\em
  Computational Mechanics}, 63:1131--1145 (2018).
\newblock

\bibitem{Bermudez2007aop}
A.~Berm\'{u}dez, L.~Hervella-Nieto, A.~Prieto, R.~Rodr\'{i}guez,
  \href{https://doi.org/10.1016/j.jcp.2006.09.018}{An optimal perfectly matched
  layer with unbounded absorbing function for time-harmonic acoustic scattering
  problems}, {\em Journal of Computational Physics}, 223:469--488 (2007).
\newblock

\bibitem{Bermudez2008aeb}
A.~Berm{\'{u}}dez, L.~Hervella-Nieto, A.~Prieto, R.~Rodr{\'{\i}}guez,
  \href{https://dx.doi.org/10.1137/060670912}{An exact bounded perfectly
  matched layer for time-harmonic scattering problems}, {\em {SIAM} Journal on
  Scientific Computing}, 30:312--338 (2008).
\newblock

\bibitem{Wu1991awr}
T.~W. Wu, A.~F. Seybert, \href{https://doi.org/10.1121/1.401901}{A weighted
  residual formulation for the {CHIEF} method in acoustics}, {\em The Journal
  of the Acoustical Society of America}, 90:1608--1614 (1991).
\newblock

\bibitem{Schenck1968iif}
H.~A. Schenck, \href{https://doi.org/10.1121/1.1911085}{Improved integral
  formulation for acoustic radiation problems}, {\em The Journal of the
  Acoustical Society of America}, 44:41--58 (1968).
\newblock

\bibitem{Nolte2014bib}
B.~Nolte, I.~Sch{\"a}fer, C.~de~Jong, L.~Gilroy, {BeTSSi II} benchmark on
  target strength simulation, in {\em Proceedings of Forum Acusticum}, 2014.

\end{thebibliography}


%
%
%
\end{document}